\documentclass[12pt,amsfonts]{article}

\usepackage{latexsym}
\usepackage{amsfonts,amssymb} 

\newtheorem{theorem}{Theorem}[section]
\newtheorem{lemma}[theorem]{Lemma}
\newtheorem{corollary}[theorem]{Corollary}
\newtheorem{proposition}[theorem]{Proposition}

\newtheorem{definition}[theorem]{Definition}
\newcommand{\bd}[1]{\begin{definition}\label{#1}\rm}
\newcommand{\ed}{\end{definition}}
\newcommand{\bt}[1]{\begin{theorem}\label{#1}}
\newcommand{\et}{\end{theorem}}
\newcommand{\bprop}[1]{\begin{proposition}\label{#1}}
\newcommand{\eprop}{\end{proposition}}
\newcommand{\bcor}[1]{\begin{corollary}\label{#1}}
\newcommand{\ecor}{\end{corollary}}

\newcommand{\T}{\textstyle}
\newcommand{\lra}{\longrightarrow}
\newcommand{\Ra}{\Longrightarrow}

\newcommand{\stack}[2]{\raisebox{-2pt} 
{\renewcommand{\arraystretch}{.01} 
\begin{tabular}{c} 
$#2$\\$\scriptscriptstyle #1$ 
\end{tabular} 
}} 

\newcommand{\wstack}[2]{\raisebox{-2pt} 
{\renewcommand{\arraystretch}{.01} 
\begin{tabular}{c} 
$\scriptscriptstyle w$\\$#2$\\$\scriptscriptstyle #1$ 
\end{tabular} 
}} 

\newcommand{\sstack}[2]{\raisebox{-2pt} 
{\renewcommand{\arraystretch}{.01} 
\begin{tabular}{c} 
$\scriptscriptstyle s$\\$#2$\\$\scriptscriptstyle #1$ 
\end{tabular} 
}}

\newcommand{\vp}{\varphi}
\newcommand{\ve}{\varepsilon}
\newcommand{\nid}{\noindent}
\newcommand{\qed}{\hfill$\Box$} 
 
\def\0{\, {\rm 0}\mskip-11mu 0} 
\def\1{\, {\rm I}\mskip-10mu 1} 
\def\g{\, \gamma\mskip-11.8mu \gamma}\,

\renewcommand{\t}[1]{\tilde{#1}} 

\newcommand{\bnu}{\mbox{\boldmath${\nu}$}} 
 
\newcommand{\sbnu}{\mbox{\scriptsize\boldmath${\nu}$}} 

\begin{document}
\title{{Mosco Type Convergence of Bilinear Forms and Weak Convergence 
of $n$-Particle Systems}} 
\par 
\author{J\"org-Uwe L\"obus 
\\ Matematiska institutionen \\ 
Link\"opings universitet \\ 
SE-581 83 Link\"oping \\ 
Sverige 
}
\date{}
\maketitle
{\footnotesize
\noindent
\begin{quote}
{\bf Abstract}
It is well known that Mosco (type) convergence is a tool in order to 
verify weak convergence of finite dimensional distributions of sequences 
of stochastic processes. In the present paper we are concerned with the 
concept of Mosco type convergence for non-symmetric stochastic processes 
and, in particular, $n$-particle systems in order to establish relative 
compactness. 
\noindent 


\end{quote}
}

\section{Introduction}
\setcounter{equation}{0}

To show weak convergence of a sequences of stochastic processes one 
has to proceed in two basic steps. These are on the one hand, proving 
relative compactness of the sequences of processes and, on the other 
hand, showing weak convergence of the finite dimensional distributions 
of the sequence of stochastic processes. 

Mosco convergence has been used in order to prove weak convergence of 
finite dimensional distributions of sequences of stochastic processes 
corresponding to symmetric Dirichlet forms, cf. for example 
\cite{CGLW12}, \cite{GKR07}, \cite{Ki06}, \cite{Ko08}, \cite{{KU97}, 
{KU96}}. 

In addition, the paper \cite{Su98} demonstrates in which way Mosco 
convergence can be used in order to verify relative compactness of 
a sequences of stochastic processes. In fact, Mosco convergence and 
additional properties or appropriate additional conditions on the 
sequences of processes provide the convergence of a certain sequence 
of associated capacities of the form ${\Bbb E}e^{-\beta\tau_n}$, 
$n\in {\Bbb N}$. Here, $\tau_n$ is a sequence of certain first exit 
times and $\beta>0$. The convergence of the sequence ${\Bbb E}e^{- 
\beta\tau_n}$, $n\in {\Bbb N}$, is then sufficient for relative 
compactness. This idea has been adapted to particle systems in 
\cite{Lo09} and \cite{Lo13} and will also be developed further in 
Section 3 of the present paper. 
\medskip 

We would also like to refer to two more motivations for this paper. 
During the last decade one may have observed an increasing interest 
in Mosco convergence relative to Dirichlet forms with changing 
reference measures or, more general, on sequences of Hilbert spaces, 
see \cite{Hi09}, \cite{Ko05}, \cite{Ko06}, \cite{Pu10}, \cite{To06}. 
Most fundamental in this sense is \cite{KS03}. In the particular 
case of sequences of $L^2$ spaces we would like to refer to 
\cite{Lo13} which has a documented history beginning 2005. 

Initiated by the two established generalizations of Dirichlet forms 
to the non-symmetric case, namely \cite{MR92} and \cite{St99}, also 
Mosco (type) convergence for non-symmetric Dirichlet forms has been 
investigated, cf. \cite{To06} and \cite{Hi98}. The paper \cite{Lo13} 
provides a framework for non-symmetric bilinear forms where neither 
the Dirichlet property nor closability is necessary. This allows to 
treat stochastic processes and particle systems without any connection 
to Dirichlet form theory. 

The present paper follows this motivation. However the framework here 
is more sophisticated. As an application, the particle system 
considered in \cite{Lo14-3} doesn't seem to be compatible with 
\cite{MR92} or \cite{St99}. Moreover the limiting initial distribution 
is no longer concentrated on one single state as in \cite{Lo13}. 
However it fits the theory presented in this article. 
\medskip 

Looking, for example, at the papers \cite{Hi98}, \cite{Hi09},
\cite{KS03}, \cite{Pu10}, \cite{To06} one can conclude that different 
classes of applications require different generalizations or 
alternations of Mosco convergence. The present paper together with 
\cite{Lo14-3} follows this attitude. 

\subsection{Basic Definitions and Technical Issues}

In order to introduce the basic setting, let ${\bnu}$ be a probability 
measure on a measurable space $(E,{\cal B})$, and let $(T_t)_{t\ge 0}$ 
be a strongly continuous contraction semigroup of linear operators on 
$L^2(E,{\bnu})$. Suppose that $(T_t)_{t\ge 0}$ is associated with a 
transition probability function $P(t,x,B)$, $t\ge 0$, $x\in E$, $B\in 
{\cal B}$, i. e., $T_tf=\int f(y)\, P(t,\cdot ,dy)$, $t\ge 0$, $f\in L^2 
(E,{\bnu})$. Assume, furthermore, that $P(t,\cdot ,E)=1$ ${\bnu}$-a.e., 
$t\ge 0$. We mention that contractivity of $(T_t)_{t\ge 0}$ on $L^2 
(E,{\bnu})$ gives $\int (P(t,x,B))^2\, \bnu(dx)\le \bnu(B)$, $B\in {\cal 
B}$, which says that $\bnu(B)=0$ implies $P(t,x,B)=0$ for $\bnu$-a.e. 
$x\in E$. 

If, as in Subsection 2.3, $(T_t)_{t\ge 0}$ is no longer contractive on 
$L^2(E,{\bnu})$ we suppose this implication to ensure well-definiteness 
of $(T_t)_{t\ge 0}$. 
\medskip

Denoting by $(A,D(A))$ the generator of $(T_t)_{t\ge 0}$ and by 
$\langle \cdot \, , \, \cdot \rangle$ the inner product in $L^2(E,{ 
\bnu})$, we introduce now the class of bilinear forms $S$ we are 
interested in. Define 
\begin{eqnarray*}
D(S):=\left\{u\in L^2(E,{\bnu}):\, \lim_{t\to 0}\left\langle 
\textstyle{\frac1t}(u-T_tu)\, , \, v\right\rangle \ \mbox{\rm exists  
for all}\ v\in L^2(E,{\bnu})\right\} 
\end{eqnarray*}
and 
\begin{eqnarray*}
S(u,v):=\lim_{t\to 0}\left\langle\textstyle{ \frac1t}(u-T_tu)\, , \, 
v\right\rangle \, , \quad u\in D(S), \ v \in L^2(E,{\bnu}). 
\end{eqnarray*}
We have $D(A)=D(S)$ according to \cite{P83}, Section 2.1 and 
\begin{eqnarray*}
S(u,v)=-\langle Au \, , \, v \rangle\, , \quad u\in D(A), \ v\in 
L^2(E,{\bnu}).
\end{eqnarray*}
In this sense we would like to understand the term {\it bilinear form}. 
However, as it is customary for Mosco (type) convergence, we also set 
$S(u,v):=\infty$ if $u\in L^2(E,{\bnu})\setminus D(S)$ and $v\in L^2 
(E,{\bnu})$. We emphasize that this definition of bilinear forms $S$ 
is adjusted to the Mosco type convergence of non-symmetric forms in 
Subsection 2.2 and, moreover, of non-positive, non-symmetric forms in 
Subsection 2.3. The latter situation appears for example in the 
application of \cite{Lo14-3}.  

Let $(G_\beta)_{\beta >0}$ be the resolvent 
associated with $S$, i. e., $G_\beta =(\beta -A)^{-1}$, $\beta >0$. 
Using contractivity of the semigroup $(T_t)_{t\ge 0}$ in $L^2(E, 
{\bnu})$ and 
\begin{eqnarray*}
\langle T_tu\, ,\, u\rangle^2\le\langle T_tu\, ,\, T_tu\rangle\langle u\, , 
\, u\rangle 
\end{eqnarray*}
one shows positivity of the form $S$, that is $S(u,u)\ge 0$ for all $u\in 
D(S)$. This observation is crucial for the whole concept of Mosco type 
convergence of sequences $S_n$ of forms on sequences of spaces $L^2(E_n, 
\bnu_n)$ to a limiting form $S$ on $L^2(E,\bnu)$ as $n\to\infty$. However, 
we also develop a framework of Mosco type convergence of sequences of 
forms when contractivity is replaced by a technical condition on $A_n'\1$, 
$n\in {\Bbb N}$, and $A'\1$ where $\1$ is the constant function taking the 
value one and the $'$ refers to the dual generator. 

The first question of interest is on the definition of such a bilinear 
form. In fact, we construct the weak generator of the semigroup $(T_t 
)_{t\ge 0}$, the second entry $v$ in $S(u,v)$ is in this sense just a 
test function. This definition of a bilinear form is far away from 
classical Dirichlet form theory which includes the notion of Mosco 
convergence. One more problem arises with the definition of the bilinear 
form, namely the appropriate notion of Mosco type convergence. We recall 
that the literature suggest several alternations of the classical one by 
U. Mosco in \cite{Mo94}, cf. \cite{Hi98}, \cite{Hi09}, and \cite{To06}. 

The results obtained in Sections 6 and 7 of \cite{Lo14-3} give rise to 
state that our approach to bilinear forms and Mosco type convergence 
is beneficial relative to the particle system we investigate there and 
similar ones. We also want to refer to a discussion on choosing the 
appropriate definition of bilinear forms relative to the mathematical 
situation, given in \cite{Lo13}, Subsections 2.1 and 2.3. 

The next problem one may come across is the existence of $A_n'\1$, 
$n\in {\Bbb N}$, and $A'\1$ in the case of non-positive bilinear forms. 
In particular, $A_n'\1$, $n\in {\Bbb N}$, and $A'\1$ should display 
properties which are useful to show Mosco type convergence, cf. Subsection 
2.3 below. This is a purely mathematical issue. In the application of 
\cite{Lo14-3}, a Fleming-Viot type particle system, this issue restricts 
the initial configurations of the particles that can be investigated in 
terms of Mosco type convergence. 
\medskip 

We conclude the introduction with a remark on the notation in the paper. 
The greek letter $\bnu$ comes always in bold. This letter is exclusively 
used to denote probability measures over spaces of probability measures, 
the states of measure valued stochastic processes. Those states of 
stochastic processes are denoted using the greek letter $\mu$, non-bold. 
 
\section{Mosco Type Convergence} 
\setcounter{equation}{0}

In this section, we are extending the framework of \cite{Lo13} in 
several ways. In Subsection 2.1, we develop the concept of convergence 
in sequences of $L^2$-spaces. This is necessary in order to establish 
the Mosco type convergence for non-symmetric forms on sequences 
of $L^2$-spaces presented in Subsection 2.2. In Subsection 2.3, we are 
then able to handle convergence of non-symmetric non-positive forms, 
the situation we have to face in the application in Section 4. 

Again, we want to refer to related research carried out in K. Kuwae, 
T. Shioya \cite{KS03} and compared with ours in \cite{Lo13}, Subsection 
3.2.  

\subsection{Analysis on Sequences of $L^2$-Spaces}

Let $\bnu_n$, $n\in {\Bbb Z}_+:=\{0,1,2,\ldots \}$, be mutually 
orthogonal probability 
measures on $(E,{\cal B})$. To ease the notation and to stress its 
special role, we will mostly use the symbol $\bnu$ instead of $\bnu_0$. 
Suppose that ${\bnu}$ is a measure with countable base on $(E,{\cal 
B})$. In addition, assume that there are mutually exclusive subsets 
$E_n$, $n\in {\Bbb Z}_+$, of $E$ such that ${\bnu}_n(E\setminus E_n) 
=0$. Let $\alpha_n$, $n\in {\Bbb Z}_+$, be a sequence of positive 
numbers with $\sum_{n=0}^\infty\alpha_n=1$. Define ${\Bbb M}:=\sum_{ 
n=0}^\infty\alpha_n{\bnu}_n$. We say that $u\in\bigcap_{n\in{\Bbb Z 
}_+}L^2(E,{\bnu}_n)$ if $u$ is an equivalence class consisting of all 
everywhere defined ${\cal B}$-measurable functions satisfying $f_1= 
f_2$ ${\Bbb M}$-a.e. if $f_1,f_2\in u$ and $\int u^2\, d{\bnu}_n < 
\infty$, $n\in {\Bbb Z}_+$. Let $\langle\cdot \, , \, \cdot\rangle_n$ 
denote the inner product in $L^2(E,{\bnu}_n)$, $n\in {\Bbb N}$, and 
let $\langle\cdot \, , \, \cdot \rangle$ denote the inner product in 
$L^2(E,{\bnu})$. 
Introduce 
\begin{eqnarray*}
{\cal D}:=\left\{\vp\in \bigcap_{\, n\in {\Bbb Z}_+}L^2(E,{\bnu}_n): 
\langle \vp\, , \, \vp\rangle_n\stack{n\to\infty}{\lra}\langle \vp 
\, , \, \vp\rangle \right\}\, . 
\end{eqnarray*}
Throughout this paper we 
suppose that there exists a linear subset ${\cal F}$ of ${\cal D}$ 
which is dense in $L^2(E,\bnu)$ and let ${\cal C}$ denote the set of 
all functions $\vp\in {\cal D}$ satisfying the following conditions: 
\begin{itemize} 
\item[(c1)] For each $\vp\in {\cal C}$, there exists a {\it 
representing sequence} $\vp_n\in {\cal F}$, $n\in {\Bbb N}$, 
such that $\vp =\vp_n$, $\bnu_n$-a.e., $n\in {\Bbb N}$. 
\item[(c2)] $\langle \vp\, , \, \psi\rangle_n\stack{n\to\infty}{\lra} 
\langle\vp\, , \, \psi\rangle$ for all $\psi\in {\cal F}$.  
\end{itemize}
Introduce 
\begin{eqnarray*}
{\cal V}:=\left\{\psi\in\bigcap_{\, n\in {\Bbb Z}_+}L^2(E,{\bnu}_n): 
\sigma\in {\cal D},\ \tau\in {\cal F}\ \mbox{\rm imply}\ \sigma\psi 
\in {\cal D},\ \tau\psi\in {\cal F}\right\}\, . 
\end{eqnarray*}
\begin{lemma}\label{Lemma2.1} 
(a) ${\cal F}\subseteq {\cal C}$. \\ 
(b) The set ${\cal C}$ is linear. \\ 
(c) The set ${\cal C}$ is dense in $L^2(E,\bnu)$. \\ 
(d) Let $\vp,\psi\in {\cal C}$. We have $\langle \vp\, , \, \psi\rangle_n 
\stack{n\to\infty}{\lra}\langle\vp\, , \, \psi\rangle$. \\ 
(e) Let $\psi\in {\cal V}$. Then, for all $\vp\in {\cal C}$, we have $\vp 
\psi\in {\cal C}$. 
\end{lemma}
Proof. (a) (c1) is trivial and (c2) follows from linearity of ${\cal 
F}$. 
\medskip

\nid
(b) Let $\ve >0$ and let $\vp,\psi\in {\cal C}$ with representing 
sequences $\vp_n\in {\cal F}$ and $\psi_n\in {\cal F}$, $n\in {\Bbb 
N}$, cf. condition (c2). Since (c1) and (c2) are obvious for $\vp +
\psi$, it remains to show that $\vp +\psi\in {\cal D}$. We can 
choose $n_0\in {\Bbb N}$ and $\t\psi\in {\cal F}$ such that, for all 
$n>n_0$, 
\begin{eqnarray}\label{2.1}
\langle\psi-\t \psi\, , \, \psi-\t \psi\rangle &<&\ve \quad 
\mbox{~since ${\cal F}\subseteq L^2(E,\bnu)$, densely,} 
\nonumber \\ 
|\langle\vp\, , \, \vp\rangle_n -\langle\vp\, , \, \vp\rangle|&<& 
\ve \quad \mbox{~since $\vp\in {\cal C}\subseteq {\cal D}$,} 
\nonumber \\
|\langle\psi\, , \, \psi\rangle_n -\langle\psi\, , \, \psi\rangle| 
&<&\ve \quad \mbox{~since $\psi\in {\cal C}\subseteq {\cal D}$,} 
\nonumber \\ 
|\langle\t \psi\, , \, \t \psi\rangle_n -\langle\t \psi\, , \, \t 
\psi\rangle|&<&\ve \quad \mbox{~since $\t \psi\in {\cal F}\subseteq  
{\cal C}\subseteq {\cal D}$,} 
\nonumber \\ 
|\langle\vp\, , \, \t \psi\rangle_n -\langle\vp\, , \, \t \psi 
\rangle|&<&\ve\, ,  \ \mbox{~cf. (c2),} 
\nonumber \\ 
|\langle\psi\, , \, \t \psi\rangle_n -\langle\psi\, , \, \t \psi 
\rangle|&<&\ve\, , \ \mbox{~cf. (c2)}\, . 
\end{eqnarray}
Therefore, we have 
\begin{eqnarray*} 
\langle \psi -\t \psi\, , \, \psi - \t \psi\rangle_n &\le& 
\langle \psi -\t \psi\, , \, \psi - \t \psi\rangle \\ 
&& \ + |\langle\psi\, , \, \psi\rangle_n -\langle\psi\, , \, \psi 
\rangle| + |\langle\t \psi\, , \, \t \psi\rangle_n -\langle\t \psi 
\, , \, \t \psi\rangle| + 2|\langle\psi\, , \, \t \psi\rangle_n 
-\langle\psi\, , \, \t \psi \rangle| \\ 
&<&5\ve
\end{eqnarray*}
which implies 
\begin{eqnarray*} 
|\langle\vp\, , \, \psi\rangle_n -\langle\vp\, , \, \t \psi 
\rangle_n|&\le&\langle\vp\, , \, \vp\rangle_n^{1/2}\langle \psi 
-\t \psi\, , \, \psi - \t \psi\rangle_n^{1/2} \\ 
&<&\sqrt{5}\sup_{n\in {\Bbb N}}\langle\vp\, , \, \vp\rangle_n^{1/2} 
\cdot\ve^{1/2} \vphantom{\int}
\end{eqnarray*}
and 
\begin{eqnarray}\label{2.2}
|\langle \vp\, , \, \psi \rangle_n - \langle \vp\, , \, \psi 
\rangle| &\le& |\langle\vp\, , \, \psi\rangle_n - \langle\vp\, , 
\, \t \psi\rangle_n| + |\langle\vp\, , \, \t \psi\rangle_n - 
\langle\vp\, , \, \t \psi\rangle| + |\langle\vp\, , \, \t \psi 
\rangle - \langle\vp\, , \, \psi\rangle| \nonumber \\ 
&<&\ve + \left(\sqrt{5}\sup_{n\in {\Bbb N}}\langle\vp\, , \, \vp 
\rangle_n^{1/2} + \langle\vp\, , \, \vp\rangle^{1/2}\right)\cdot 
\ve^{1/2}\, . 
\end{eqnarray}
It follows now from (\ref{2.1}) and (\ref{2.2}) that 
\begin{eqnarray*}
|\langle\vp+\psi\, , \, \vp+\psi\rangle_n-\langle\vp+\psi\, , \, 
\vp+\psi\rangle| \\ 
&&\hspace{-3cm}\le|\langle\vp\, , \, \vp\rangle_n-\langle\vp\, , \, 
\vp\rangle|+|\langle\psi\, , \, \psi\rangle_n-\langle\psi\, , \, 
\psi\rangle|+2|\langle \vp\, , \, \psi \rangle_n - \langle \vp\, , 
\, \psi\rangle| \\ 
&&\hspace{-3cm}< 4\ve + \left(2\sqrt{5}\sup_{n\in {\Bbb N}}\langle 
\vp\, , \, \vp\rangle_n^{1/2} + 2\langle\vp\, , \, \vp\rangle^{1/2} 
\right)\cdot \ve^{1/2}\, . 
\end{eqnarray*}
(c) This is a consequence of (a) and the fact that ${\cal F}$ is 
dense in $L^2(E,\bnu)$. 
\medskip 

\nid 
(d) This follows from (\ref{2.2}). 
\medskip 

\nid 
(e) Let $\vp\in {\cal C}$ and $\psi\in{\cal V}$. By hypothesis, we have 
$\vp\psi\in {\cal D}$. Let $\vp_n\in {\cal F}$, $n\in {\Bbb N}$, be the 
representing sequence of $\vp$, cf. (c1). Then $\vp_n\psi\in {\cal F}$, $n 
\in {\Bbb N}$, by hypothesis. In other words, $\vp_n\psi\in {\cal F}$, $n 
\in {\Bbb N}$, is the representing sequence of $\vp\psi$, i. e., we have (c1) 
for $\vp\psi$. Furthermore, for all $\rho\in {\cal F}$, we have $\psi\rho 
\in {\cal F}$ by hypothesis and, by (c2), $\langle \vp\, , \, \psi\rho 
\rangle_n\stack{n\to\infty}{\lra}\langle\vp\, , \, \psi\rho\rangle$. Thus, 
$\langle\vp\psi\, , \, \rho\rangle_n\stack{n\to\infty}{\lra}\langle\vp\psi 
\, , \, \rho\rangle$. Thus, we have (c2) for $\vp\psi$. 
\qed 
\medskip 

\nid
\begin{definition}\label{Definition2.2} {\rm 
(a) A sequence $\vp_n\in {\cal C}$, $n\in {\Bbb N}$, is said to 
be {\it w-convergent} to $\vp\in L^2(E,{\bnu})$ as $n\to\infty$ 
(in symbols $\vp_n\wstack{n\to\infty}{\lra} \vp$) if 
\begin{itemize} 
\item[(i)] $\langle\vp_n\, , \, \psi\rangle_n\stack{n\to 
\infty}{\lra}\langle\vp\, , \, \psi\rangle$ for all $\psi\in 
{\cal C}$.  
\end{itemize} 
(b) A sequence $\psi_n\in {\cal C}$, $n\in {\Bbb N}$, is said to 
be {\it s-convergent} to $\psi\in L^2(E,{\bnu})$ as $n\to 
\infty$ (in symbols $\psi_n\sstack{n\to\infty}{\lra} \psi$) if 
\begin{itemize} 
\item[(i)] $\psi_n$ $w$-converges to $\psi$ as $n\to\infty$ and 
\item[(ii)] $\langle \psi_n\, , \, \psi_n\rangle_n \stack{n\to 
\infty}{\lra}\langle \psi\, , \, \psi\rangle$.
\end{itemize} 
(c) Speaking of $w$-convergence or $s$-convergence of subsequences 
$\vp_{n_k}\in {\cal C}$ or $\psi_{n_k}\in {\cal C}$, respectively, 
will mean that in (a) or (b) the index $n\in {\Bbb N}$ is replaced 
with $n_k\in {\Bbb N}$. }
\end{definition}
{\bf Remarks} (1) Let $\psi_n\in {\cal C}$, $n\in {\Bbb N}$, 
$w$-converge to $\psi\in {\cal C}$. According to the definition of 
${\cal C}$, that $\psi_n$ $s$-converges to $\psi\in {\cal C}$ is 
equivalent to 
\begin{eqnarray*}
\langle\psi_n-\psi \, , \, \psi_n-\psi\rangle_n &=&\langle\psi_n 
\, , \, \psi_n\rangle_n -2\langle\psi_n\, , \, \psi\rangle_n + 
\langle\psi \, , \, \psi\rangle_n \\ 
&\stack{n\to\infty}{\lra}& 0\, . 
\end{eqnarray*}
(2) Let $\psi\in {\cal C}$. It follows from Lemma \ref{Lemma2.1} (d), 
that $\psi_n:=\psi$, $n\in {\Bbb N}$, $s$-converges to $\psi$. 
\begin{proposition}\label{Proposition2.3} 
(a) Let $\vp_n\in {\cal C}$, $n\in {\Bbb N}$, be a sequence 
$w$-convergent to $\vp\in L^2(E,{\bnu})$ as $n\to\infty$. Then 
$\langle \vp_n\, , \, \vp_n\rangle_n$, $n\in {\Bbb N}$, is bounded. \\ 
(b) Let $\vp_n\in {\cal C}$, $n\in {\Bbb N}$, be a sequence 
such that $\langle \vp_n\, , \, \vp_n\rangle_n$ is bounded. Then 
there exists a subsequence $\vp_{n_k}\in {\cal C}$, $k\in {\Bbb N}$, 
$w$-convergent to some $\vp\in L^2(E,{\bnu})$ as $k\to\infty$. \\  
(c) Let $\vp_n\in {\cal C}$, $n\in {\Bbb N}$, be a sequence 
$w$-convergent to $\vp\in L^2(E,{\bnu})$ and let $\psi_n\in {\cal C}$, $n 
\in {\Bbb N}$, be a sequence that $s$-converges to $\psi\in L^2(E,{\bnu})$ 
as $n\to\infty$. Then $\langle\vp_n\, , \, \psi_n\rangle_n \stack{n\to 
\infty}{\lra}\langle\vp\, , \, \psi\rangle$. \\ 
(d) Let $\vp_n\in {\cal C}$, $n\in {\Bbb N}$, be a sequence $w$-convergent 
to $\vp\in {\cal C}$ as $n\to\infty$ and let $\psi_n\in {\cal C}$, $n\in 
{\Bbb N}$, be a sequence $w$-convergent to $\psi\in {\cal C}$ as $n\to 
\infty$. Suppose 
\begin{eqnarray}\label{2.3}
\langle \vp_n -\vp\, , \, \psi_n-\psi \rangle_n \ge 0\, , \quad 
n\in {\Bbb N}. 
\end{eqnarray} 
Then 
\begin{eqnarray*}
\liminf_{n\to\infty}\langle\vp_n\, , \, \psi_n\rangle_n\ge \langle\vp 
\, , \, \psi\rangle\, . 
\end{eqnarray*} 
(e) Let $\vp_n\in {\cal C}$, $n\in {\Bbb N}$, be a sequence $w$-convergent 
to $\vp\in L^2(E,{\bnu})$. Then we have 
\begin{eqnarray*}
\liminf_{n\to\infty}\langle\vp_n \, ,\, \vp_n\rangle_n\ge \langle 
\vp\, , \, \vp\rangle\, . 
\end{eqnarray*}

In the following, let $\psi_n\in{\cal V}$, $n\in {\Bbb N}$, let $\psi\in 
L^\infty (E,\bnu)$, and assume $\psi_n\rho\sstack{n\to\infty}{\lra}\psi 
\rho$ for all $\rho\in {\cal C}$. 
\smallskip 

\nid
(f) Let ${\cal C}\ni\vp_n\wstack{n\to\infty}{\lra}\vp\in L^2(E,\bnu)$. 
Then $\vp_n\psi_n\wstack{n\to\infty}{\lra}\vp\psi$. \\ 
(g) Let ${\cal C}\ni\vp_n\sstack{n\to\infty}{\lra}\vp\in L^2(E,\bnu)$. 
Then $\vp_n\psi_n\sstack{n\to\infty}{\lra}\vp\psi$. 
\end{proposition}
Proof. For the proofs of (a) through (c), we will refer to \cite{Lo13}. 
\medskip 

\nid
(a) and (c) The situation in Section 3 of \cite{Lo13} is compatible 
with the setting here. In particular, there it is assumed that $E$ is 
a metric space and Definition 3.1 yields a linear set ${\cal C}$. 
Replacing in the statement of \cite{Lo13}, Lemma 3.2 (b), $C_b(E)$ by 
${\cal F}$ and using the above properties of ${\cal F}$ the assertion 
of \cite{Lo13}, Lemma 3.2 (b), becomes obvious. The proofs of 
Proposition 3.3 (a) and (b) can now be followed word for word with 
${\cal F}$ instead of $C_b(E)$. In this way, we have verified (a) and 
(c) of the present proposition. 
\medskip 

\nid
(b) This has been demonstrated in \cite{Lo13}, proof of Proposition 
2.3 (a). 
\medskip 

\nid 
(d) Because of $\vp_n\wstack{n\to\infty}{\lra}\vp$, $\psi_n\wstack 
{n\to\infty}{\lra}\psi$, and $\vp ,\psi \in {\cal C}$, we have $\langle 
\vp_n\, , \, \psi \rangle_n\stack{n\to\infty}{\lra}\langle\vp\, , \, 
\psi\rangle$, $\langle\vp\, , \, \psi_n \rangle_n\stack{n\to\infty} 
{\lra}\langle\vp\, , \, \psi\rangle$, and $\langle\vp\, , \, \psi 
\rangle_n\stack{n\to\infty}{\lra}\langle\vp\, , \, \psi\rangle$. The 
lemma is now a consequence of hypothesis (\ref{2.3}) and 
\begin{eqnarray*}
\langle\vp_n\, , \, \psi_n\rangle_n&=&\langle\vp_n\, , \, \psi\rangle_n 
+\langle\vp\, , \, \psi_n\rangle_n-\langle\vp\, , \, \psi\rangle_n  
+\langle \vp_n -\vp\, , \, \psi_n-\psi \rangle_n\, . 
\end{eqnarray*} 
(e) Let $\ve>0$ and $\t\vp\in {\cal C}$ such that $\langle\t\vp- 
\vp \, ,\, \t\vp-\vp\rangle<\ve$, cf. Lemma \ref{Lemma2.1} (c). 
Since $\vp_n$, $n\in {\Bbb N}$, $w$-converges to $\vp\in L^2(E, 
{\bnu})$, it follows that  
\begin{eqnarray*}
&&\hspace{-.5cm}\langle\vp_n\, ,\, \vp_n\rangle_n-\langle\vp_n-\t 
\vp\, ,\, \vp_n-\t\vp\rangle_n=\langle\vp_n\, ,\, \t\vp\rangle_n 
+\langle\t\vp\, ,\, \vp_n\rangle_n-\langle\t\vp\, ,\, \t\vp\rangle_n 
 \\ 
&&\stack{n\to\infty}{\lra}\langle\vp\, ,\, \t\vp\rangle+\langle\t 
\vp\, ,\, \vp\rangle-\langle\t\vp\, ,\, \t\vp\rangle \\ 
&&\hspace{.5cm}=\langle\vp\, ,\, \vp\rangle-\langle\vp-\t\vp\, ,\, 
\vp-\t\vp\rangle\, .  
\end{eqnarray*}
We have $\langle\vp-\t\vp\, ,\, \vp-\t\vp\rangle\le\ve$ from which 
we get $\liminf_{n\to\infty}\langle\vp_n\, ,\, \vp_n\rangle_n\ge 
\langle \vp\, , \, \vp\rangle$. 
\medskip 

\nid 
(f) Let $\rho\in {\cal C}$. According to Lemma \ref{Lemma2.1} (e) and 
hypotheses, we have $\vp_n\psi_n\in {\cal C}$, $n\in {\Bbb N}$, and ${\cal 
C}\ni\psi_n\rho\sstack{n\to\infty}{\lra}\psi\rho\in L^2(E,\bnu)$. From part 
(c), we 
obtain  
\begin{eqnarray*}
\langle\vp_n\psi_n\, , \, \rho\rangle_n=\langle\vp_n\, , \, \psi_n\rho 
\rangle_n\stack{n\to\infty}{\lra}\langle\vp\, , \, \psi\rho\rangle=\langle 
\vp\psi\, , \, \rho\rangle\, . 
\end{eqnarray*} 
(g) It remains to note that by parts (c) and (f), we have $\vp_n\psi_n^2 
\wstack{n\to\infty}{\lra}\vp\psi^2$ which yields by part (c)
\begin{eqnarray*}
\langle\vp_n\psi_n\, , \, \vp_n\psi_n\rangle&&\hspace{-.5cm}=\langle\vp 
\psi_n\, , \, \vp_n\psi_n\rangle_n+\langle\vp_n-\vp\, , \, \vp_n\psi_n^2 
\rangle_n \\ 
&&\hspace{-1cm}\stack{n\to\infty}{\lra}\langle\vp\psi\, , \, \vp\psi 
\rangle\, . 
\end{eqnarray*} 
\qed 

\subsection{Mosco Type Convergence of Non-Symmetric Forms on Sequences 
of $L^2$-Spaces}
 
For every $n\in {\Bbb N}$, let $(T_{n,t})_{t\ge 0}$ be a strongly 
continuous contraction semigroup in $L^2(E,\bnu_n)$ and let $(T_t 
)_{t\ge 0}$ be a strongly continuous contraction semigroup in $L^2 
(E,\bnu)$, all in the sense of Section 1. Denote by $S_n$, $A_n$, 
$(G_{n,\beta})_{\beta >0}$ the bilinear form in the sense of 
Section 1, the generator, and the family of resolvents associated 
with $(T_{n,t})_{t\ge 0}$, $n\in {\Bbb N}$. Similarly, let $S$, 
$A$, and $(G_\beta)_{\beta >0}$ the bilinear form, the generator, 
and the family of resolvents associated with $(T_t)_{t\ge 0}$. For 
$\beta >0$, $n\in {\Bbb N}$, $\vp_n\in D(S_n)$, $\psi_n\in L^2(E, 
{\bnu}_n)$, set $S_{n,\beta}(\vp_n,\psi_n):=\beta\langle\vp_n\, , 
\, \psi_n \rangle_n + S_n(\vp_n,\psi_n)$, and for $\vp\in D(S)$, 
$\psi\in L^2(E,{\bnu})$, define $S_{\beta}(\vp,\psi):=\beta\langle 
\vp\, ,\, \psi\rangle + S(\vp,\psi)$. 
\begin{definition}\label{Definition2.4} 
{\rm We say that $S_n$, $n\in {\Bbb N}$, {\it pre-converges} to $S$ if 
\begin{itemize} 
\item[(i)] For every $\vp\in L^2(E,{\bnu})$ and every subsequence 
$\vp_{n_k}\in D(S_{n_k})\cap {\cal C}$, $k\in {\Bbb N}$, $w$-converging 
to $\vp$ such that $\sup_{k\in {\Bbb N}}\left\langle A_{n_k}\vp_{n_k}\, 
,\, A_{n_k}\vp_{n_k}\right\rangle_{n_k}<\infty$, we have 
\begin{eqnarray*}
S(\vp,\vp )\le \liminf_{k\to\infty}S_{n_k}(\vp_{n_k},\vp_{n_k} )\, . 
\end{eqnarray*}
\item[(ii)] For every $\psi\in D(S)$, there exists a sequence $\psi_n 
\in D(S_n)\cap {\cal C}$, $n\in {\Bbb N}$, $s$-converging to $\psi$ 
such that $\sup_{n\in {\Bbb N}}\left\langle A_n\psi_n\, ,\, A_n\psi_n 
\right\rangle_n<\infty$ and 
\begin{eqnarray*}
\limsup_{n\to\infty} S_n(\psi_n ,\psi_n )\le S(\psi ,\psi )\, .  
\end{eqnarray*}
\end{itemize} 
}
\end{definition}
\begin{lemma}\label{Lemma2.5} 
Let $S_n$, $n\in {\Bbb N}$, be a sequence of bilinear forms 
pre-convergent to $S$. In addition, let $w_n\in D(S_n)\cap {\cal C}$, 
$n\in {\Bbb N}$, be a sequence $w$-converging to  some $w\in L^2(E, 
{\bnu})$ and let $v_n\in D(S_n)\cap {\cal C}$, $n\in {\Bbb N}$, be 
a sequence $s$-converging to $v\in D (S)$ in the sense of condition 
(ii) of Definition \ref{Definition2.4}. Suppose $\, \sup_{n\in {\Bbb 
N}}\langle A_nw_n\, ,\, A_nw_n\rangle_n<\infty$. \\ 
(a) The limit $\lim_{n\to\infty}(S_n(v_n,w_n)+S_n(w_n,v_n))$ exists 
and we have 
\begin{eqnarray*}
\lim_{n\to\infty}(S_n(v_n,w_n)+S_n(w_n,v_n))=S(v,w)+S(w,v)\, . 
\end{eqnarray*}
(b) For $\beta >0$ and $\, \limsup_{n\to\infty}S_{n,\beta}(w_n,w_n)< 
\infty$ then 
\begin{eqnarray*}
\lim_{n\to\infty}(S_{n,\beta}(v_n,w_n)+S_{n,\beta}(w_n,v_n))=S_\beta 
(v,w)+S_\beta (w,v)\, . 
\end{eqnarray*}
The lemma holds also for subsequences $n_k$, $k\in {\Bbb N}$, of indices. 
\end{lemma}
Proof. See \cite{Lo13}, proof of Lemma 2.5. However note the difference 
in the definition of pre-convergence, cf. condition (i) in Definition 
\ref{Definition2.4}. Note also that $\, \sup_{n\in {\Bbb N}}\langle A_n 
w_n\, ,\, A_nw_n\rangle_n<\infty$ implies by Proposition 
\ref{Proposition2.3} (a) $\limsup_{n\to\infty}S_n(w_n,w_n)<\infty$, cf. 
the corresponding assumptions of Lemma 2.5 in \cite{Lo13}. 

Also in the present exposition, we observe that the pre-convergence 
of $S_n$, $n\in {\Bbb N}$, to $S$ implies the pre-convergence of 
$S_{n,\beta}$, $n\in {\Bbb N}$, to $S_\beta$. Property (i) of Definition 
\ref{Definition2.4} for $S_{n,\beta}$, $n\in {\Bbb N}$, and $S_\beta$ 
follows from Proposition \ref{Proposition2.3} (a) and (e). The sequence 
$\psi_n$, $n\in {\Bbb N}$, in property (ii) of Definition 
\ref{Definition2.4} for $S_{n,\beta}$, $n\in {\Bbb N}$, and $S_\beta$  
is the same as that in property (ii) for $S_n$, $n\in {\Bbb N}$, and $S$. 
\qed
\medskip 

Let us introduce the following condition. 
\begin{itemize}
\item[(c3)] 
\begin{itemize}
\item[$(i)$] ${\cal G}:=\{G_\beta g:g\in {\cal C},\ \beta >0\}\subseteq 
{\cal C}$ in the sense that for every $g \in {\cal C}$ and $\beta>0$, 
there is a $u\in {\cal C}$ with $G_\beta g=u$ ${\bnu}$-a.e.
\item[$(ii)$] ${\cal G}_n:=\{G_{n, \beta}g:g\in {\cal C},\ \beta >0\} 
\subseteq {\cal C}$, $n\in {\Bbb N}$, in the sense that for every $g 
\in {\cal C}$, $\beta>0$, and every $n\in {\Bbb N}$, there exists a $v 
\in {\cal C}$ such that $G_{n,\beta}g=v$ ${\bnu}_n$-a.e. 
\item[$(iii)$] ${\cal G}':=\{G'_\beta g:g\in {\cal C},\ \beta >0\}\subseteq 
{\cal C}$. 
\item[$(iv)$] ${\cal G}'_n:=\{G'_{n, \beta}g:g\in {\cal C},\ \beta >0\} 
\subseteq {\cal C}$, $n\in {\Bbb N}$. 
\end{itemize}
\end{itemize}
{\bf Remark} (3) Imposing condition (i) in Definition \ref{Definition2.4} 
on $S_n$, $n\in {\Bbb N}$, and $S$, we implicitly require that $\vp\in D 
(S)$. On the one hand we suppose that $\vp_{n_k}\in D(S_{n_k})\cap {\cal 
C}$, $k\in {\Bbb N}$, $w$-converges to $\vp$. This implies $\sup_{k\in 
{\Bbb N}}\left\langle\vp_{n_k}\, ,\, \vp_{n_k}\right\rangle_{n_k}<\infty$ 
by Proposition \ref{Proposition2.3} (a). On the other hand, we require 
$\sup_{k\in {\Bbb N}}\left\langle A_{n_k}\vp_{n_k}\, ,\, A_{n_k}\vp_{n_k} 
\right\rangle_{n_k}<\infty$. Thus 
\begin{eqnarray*}
S(\vp,\vp )\le\liminf_{k\to\infty}S_{n_k}(\vp_{n_k},\vp_{n_k} )\le 
\limsup_{k\to\infty}\left\langle -A_{n_k}\vp_{n_k}\, ,\, \vp_{n_k} 
\right\rangle_{n_k}<\infty
\end{eqnarray*}
which, by definition, says that $\vp\in D(S)$. 

Conversely, in order to verify condition (i) in Definition 
\ref{Definition2.4}, it makes sense to show that $D(S_{n_k})\cap {\cal C} 
\ni\vp_{n_k}\wstack{n\to\infty}{\lra}\vp\in L^2(E,\bnu)$ and $\sup_{k\in 
{\Bbb N}}\left\langle A_{n_k}\vp_{n_k}\, ,\, A_{n_k}\vp_{n_k}\right 
\rangle_{n_k}<\infty$ imply $\vp\in D(S)=D(A)$. This is how we proceed 
in the proof of \cite{Lo14-3}, Proposition 6.1, Step 3. 
\begin{lemma}\label{Lemma2.6} 
Let $S_n$, $n\in {\Bbb N}$, be a sequence of bilinear forms 
pre-convergent to $S$. Furthermore, let $\beta >0$ and let $u_n\in 
D(S_n)\cap {\cal C}$ such that $A_nu_n\in {\cal C}$, $n\in {\Bbb N}$, be a 
$w$-convergent sequence with $\, \sup_{n\in {\Bbb N}}\langle A_nu_n\, ,\, 
A_nu_n\rangle_n<\infty$. Let $u\in D(S)$. Introduce the following 
conditions. 
\begin{itemize}
\item[{\rm (iii)}] Let $u_n$, $n\in {\Bbb N}$, and $u$ as above.  
\begin{eqnarray}\label{2.4} 
\lim_{n\to\infty}S_{n,\beta}(u_n,\psi_n)=S_\beta (u,\psi)
\end{eqnarray}
for all $\psi\in {\cal C}$ and all sequences $\psi_n\in {\cal C}$, $n\in 
{\Bbb N}$, $s$-convergent to $\psi$ yields 
\begin{eqnarray}\label{2.5}
\lim_{n\to\infty}S_{n,\beta}(\psi_n,u_n)=S_\beta (\psi ,u)  
\end{eqnarray}
for all $\psi\in D(S)$ and all sequences $\psi_n\in D(S_n)\cap {\cal 
C}$, $n\in {\Bbb N}$, $s$-convergent to $\psi$ in the sense of condition 
(ii) in Definition \ref{Definition2.4}. 
\item[{\rm (iv)}] Let $u_n$, $n\in {\Bbb N}$, and $u$ as above. If 
$\beta u_n-A_nu_n\wstack{n\to\infty}{\lra}\beta u-Au$ then $u_n\wstack 
{n\to\infty}{\lra}u$. 
\end{itemize} 
Then (iv) implies (iii). 
\end{lemma}
\medskip

\nid
{\bf Remark} (4) In (iii), we require that (\ref{2.4}) holds for all $\psi 
\in {\cal C}$ and all sequences $\psi_n\in {\cal C}$, $n\in {\Bbb N}$, 
$s$-convergent to $\psi$. It is equivalent to replace in this sentence 
$\psi\in {\cal C}$ by $\psi\in {\cal C}'$ for any dense subset ${\cal C}'$ 
of $L^2(E,\bnu)$. For this, recall also Proposition \ref{Proposition2.3} (a). 
\medskip

\nid
Proof. Let us assume (iv) and (\ref{2.4}). We verify (\ref{2.5}). 
For this, let us specify for a moment $\psi_n:=\psi$, $n\in {\Bbb N}$. 
It follows then from Remark (2) and (\ref{2.4}) that $\beta u_n-A_n 
u_n\wstack{n\to\infty}{\lra}\beta u-Au$. Condition (iv) implies now 
that $u_n\wstack{n\to\infty}{\lra}u$. From Lemma \ref{Lemma2.5} and 
(\ref{2.4}), we obtain 
\begin{eqnarray*}
\lim_{n\to\infty}S_{n,\beta}(\psi_n,u_n)&=&\lim_{n\to\infty}\left\{ 
S_{n,\beta}(\psi_n,u_n)+S_{n,\beta}(u_n,\psi_n)\right\}-\lim_{n\to 
\infty}S_{n,\beta}(u_n,\psi_n) \\ 
&=&\left\{ S_{\beta}(\psi,u)+S_{\beta}(u,\psi)\right\}-S_{\beta}(u, 
\psi) \\ 
&=&S_{\beta}(\psi,u) \vphantom{\sum}
\end{eqnarray*}
for all $\psi\in D(S)$ and all sequences $\psi_n\in D(S_n)\cap {\cal 
C}$, $n\in {\Bbb N}$, $s$-convergent to $\psi$ satisfying condition 
(ii) in Definition \ref{Definition2.4}. 
\qed
\medskip 

\nid
{\bf Definition 2.4 continued } Let $S_n$, $n\in {\Bbb N}$, be a 
sequence of bilinear forms pre-convergent to $S$. If, in addition, 
condition (iii) in Lemma \ref{Lemma2.6} is satisfied, then we say 
that $S_n$, $n\in {\Bbb N}$, {\it converges} to $S$. 
\begin{theorem}\label{Theorem2.7}
Let $\beta>0$, suppose that conditions (c1)-(c3) are satisfied, and 
assume that $S_n$, $n\in {\Bbb N}$, converges to $S$ in the sense 
of Definition \ref{Definition2.4}. \\ 
(a) For all $f\in L^2(E,\bnu)$ and all sequences $f_n\in {\cal C}$ 
$w$-converging to $f$, $G_{n,\beta}f_n$ $w$-converges to $G_\beta f$ 
and $G_{n,\beta}'f_n$ $w$-converges to $G_\beta'f$ as $n\to\infty$. \\ 
(b) For all $g\in L^2(E,\bnu)$ and all sequences $g_n\in {\cal C}$ 
$s$-converging to $g$, $G_{n,\beta}g_n$ $s$-converges to $G_\beta g$ 
and $G_{n,\beta}'g_n$ $s$-converges to $G_\beta'g$ as $n\to\infty$. 
\end{theorem}
Proof. In Step 1 below, we will show that for all $g\in {\cal C}$, 
$G_{n,\beta}g$ $s$-converges to $G_\beta g$. In Step 2 we will use 
the ideas of Step 1 to prove that for all $f\in L^2(E,\bnu)$ and all 
sequences $f_n\in {\cal C}$ $w$-converging to $f$, $G_{n,\beta}f_n$ 
$w$-converges to $G_\beta f$. In Step 3, we will demonstrate that 
for all $g\in L^2(E,\bnu)$ and all sequences $g_n\in {\cal C}$ 
$s$-converging to $g$, $G_{n,\beta}g_n$ $s$-converges to $G_\beta g$. 
Finally, Step 4 will be devoted to the verification of the second 
part of (a) as a consequence of the first part of (b). A straight 
forward conclusion will then be the second part of (b).
\medskip 

\nid 
{\it Step 1 } Fix $g\in {\cal C}$ and $\beta >0$. Set $u_n:=G_{n, 
\beta}g$. Because of $\langle u_n\, ,\, u_n\rangle_n^{1/2}\le\frac 
{1}{\beta}\sup_{n'\in {\Bbb N}}\langle g\, , \, g\rangle^{1/2}_{n'} 
<\infty$, $n\in {\Bbb N}$, and Proposition \ref{Proposition2.3} (b), 
there exists a subsequence $u_{n_k}$, $k\in {\Bbb N}$, $w$-converging 
to some $\t u\in L^2(E,\bnu)$. For this, recall also condition 
(c3). Because of Remark (3), we even may conclude $\t u\in D(S)$ since 
\begin{eqnarray*}
\sup_{k\in {\Bbb N}}\left\langle A_{n_k}u_{n_k}\, ,\, A_{n_k}u_{n_k} 
\right\rangle_{n_k}&&\hspace{-.5cm}=\sup_{k\in {\Bbb N}}\left\langle 
\beta u_{n_k}-g\, ,\, \beta u_{n_k}-g\right\rangle_{n_k}\le 4\sup_{k 
\in {\Bbb N}}\langle g\, ,\, g\rangle_{n_k}<\infty\, .  
\end{eqnarray*}

Set $u:=G_{\beta}g$ and let $\psi\in D(S)$. We have $\lim_{n\to\infty} 
S_{n,\beta}(G_{n,\beta}g,\psi_n)=\lim_{n\to\infty}\langle g\, ,\, 
\psi_n\rangle_n =\langle g\, ,\, \psi\rangle=S_{\beta}(G_{\beta}g,\psi 
)$ for all sequences $\psi_n\in{\cal C}$ $s$-converging to $\psi$. 
Thus, condition (iii) of Definition \ref{Definition2.4}, Remark (4), 
and Lemma \ref{Lemma2.5} imply that 
\begin{eqnarray}\label{2.6}
S_\beta (\psi,u)+S_\beta(u,\psi)&=&\lim_{k\to\infty}\left\{S_{n_k, 
\beta}(\psi_{n_k},G_{n_k,\beta}g)+S_{n_k,\beta}(G_{n_k,\beta}g, 
\psi_{n_k})\right\}\nonumber \\ 
&=&S_\beta (\psi,\t u)+S_\beta(\t u,\psi) 
\end{eqnarray}
for all $\psi\in D(S)$ and all sequences $\psi_n\in D(S_n)\cap {\cal 
C}$ $s$-converging to $\psi$ in the sense of condition (ii) of 
Definition \ref{Definition2.4}. Note that, in order to use Lemma 
\ref{Lemma2.5}, we verify $\sup_{n\in {\Bbb N}}\langle A_nu_n\, ,\, 
A_nu_n\rangle_n<\infty$ as above. 

Applying (\ref{2.6}) to both, $\psi=u$ and $\psi=\t u$, we conclude 
$S_\beta (u-\t u,u-\t u)=0$ and thus $u=\t u$. In other words, $G_{n, 
\beta}g\wstack{n\to\infty}{\lra}G_\beta g$, independent of the above 
chosen subsequence $n_k$, $k\in {\Bbb N}$. $G_{n,\beta} g\sstack{n\to 
\infty}{\lra}G_\beta g$ is now a consequence of Proposition 
\ref{Proposition2.3} (e), condition (c3) $(ii)$, and Definition 
\ref{Definition2.4} (i), which imply   
\begin{eqnarray*} 
\beta\langle G_\beta g\, , \, G_\beta g\rangle +S(G_\beta g,G_\beta g) 
&\le&\liminf_{n\to\infty}\left\{\beta\langle G_{n,\beta}g\, , \, G_{n, 
\beta}g\rangle_n +S_n(G_{n,\beta} g,G_{n,\beta} g)\right\} \\ 
&=&\liminf_{n\to\infty}\langle g\, , \, G_{n,\beta} g\rangle_n \\ 
&=&\langle g\, , \, G_\beta g\rangle \\ 
&=&\beta\langle G_\beta g\, , \, G_\beta g\rangle +S(G_\beta g,G_\beta 
g)
\end{eqnarray*}
and thus $\lim_{n\to\infty}\langle G_{n,\beta}g\, ,\, G_{n,\beta}g 
\rangle_n=\langle G_\beta g\, ,\, G_\beta g\rangle$. 
\medskip

\nid 
{\it Step 2 } Let $f\in L^2(E,\bnu)$, $f_n\in {\cal C}$, $n\in{\Bbb N}$,
be a sequence $w$-converging to $f$. Set $u_n:=G_{n,\beta}f_n$. By 
Proposition \ref{Proposition2.3} (a) we have $\langle u_n\, ,\, u_n 
\rangle_n^{1/2}\le\frac{1}{\beta}\sup_{n'\in{\Bbb N}}\langle f_{n'}\, ,\, 
f_{n'}\rangle^{1/2}_{n'}<\infty$, $n\in {\Bbb N}$. Because of Proposition  
\ref{Proposition2.3} (b) and Remark (3) there exists a subsequence $u_{ 
n_k}$ $w$-converging to some $\t u\in D(S)$ as $k\to\infty$. Here we have 
used %
\begin{eqnarray*}
\sup_{k\in {\Bbb N}}\left\langle A_{n_k}u_{n_k}\, ,\, A_{n_k}u_{n_k} 
\right\rangle_{n_k}&&\hspace{-.5cm}=\sup_{k\in {\Bbb N}}\left\langle 
\beta u_{n_k}-f_{n_k}\, ,\, \beta u_{n_k}-f_{n_k}\right\rangle_{n_k}\le 
4\sup_{k\in {\Bbb N}}\langle f_{n_k}\, ,\, f_{n_k}\rangle_{n_k}<\infty\, , 
\end{eqnarray*}
the latter by Proposition \ref{Proposition2.3} (a). 

Set $u:=G_{\beta}f$ and let $\psi\in D(S)$. As in Step 1, we have $\lim_{n 
\to\infty}S_{n,\beta}(G_{n,\beta}f_n,\psi_n)=\lim_{n\to\infty}\langle f_n 
\, ,\, \psi_n\rangle_n=\langle f\, ,\, \psi\rangle =S_{\beta}(G_{\beta}f, 
\psi)$ for all sequences $\psi_n\in{\cal C}$ $s$-converging to $\psi$. Thus, 
condition (iii) of Definition \ref{Definition2.4}, Remark (4), and Lemma 
\ref{Lemma2.5} yield as in Step 1 $S_\beta (\psi,u)+S_\beta(u,\psi)= 
S_\beta(\psi,\t u)+S_\beta(\t u,\psi)$ for all $\psi\in D(S)$. Note again 
that, in order to use Lemma \ref{Lemma2.5}, we have $\sup_{k\in {\Bbb N}} 
\left\langle A_{n_k}u_{n_k}\, ,\, A_{n_k}u_{n_k}\right\rangle_{n_k}<\infty$. 
Again we may conclude $u=\t u$. We have thus verified $G_{n,\beta}f_n\wstack
{n\to\infty}{\lra}G_\beta f$, independent of the above chosen subsequence 
$n_k$, $k\in{\Bbb N}$. 
\medskip

\nid 
{\it Step 3 } Now, let $g\in L^2(E,\bnu)$ and $g_n\in {\cal C}$, $n\in 
{\Bbb N}$, such that $g_n\sstack{n\to\infty}{\lra}g$. Let $\vp\in {\cal 
C}$, $\ve >0$, and choose $\t g\in {\cal C}$ with $\langle g-\t g\, ,\, 
g-\t g\rangle^{1/2}<\ve$. We have 
\begin{eqnarray*}
&&\hspace{-1cm}\left|\left(\langle G_{n,\beta}g_n\, , \, \vp\rangle_n- 
\langle G_\beta g\, , \, \vp\rangle\right)-\left(\langle G_{n,\beta}\t g 
\, , \, \vp\rangle_n-\langle G_\beta\t g\, , \, \vp\rangle\right)\right| \\ 
&&\le \left|\langle G_{n,\beta}(g_n-\t g)\, , \, \vp\rangle_n\right|+\left| 
\langle G_\beta (g-\t g)\, , \, \vp\rangle\right| \vphantom{\sum}\\ 
&&\le {\textstyle\frac1\beta}\langle\vp\, , \, \vp\rangle^{1/2}_n\left( 
\langle g_n\, , \, g_n\rangle_n-2\langle g_n\, , \, \t g\rangle_n +\langle 
\t g\, , \, \t g\rangle_n\right)^{1/2}+{\textstyle \frac{\ve}{\beta}}\langle 
\vp\, , \, \vp\rangle^{1/2} \\ 
&&\hspace{-.5cm}\stack{n\to\infty}{\lra}{\textstyle\frac1\beta}\langle\vp 
\, , \, \vp\rangle^{1/2}\langle g-\t g\, , \, g-\t g\rangle^{1/2} 
+{\textstyle\frac{\ve}{\beta}}\langle\vp\, , \, \vp\rangle^{1/2} 
={\textstyle\frac{2\ve}{\beta}}\langle\vp\, , \, \vp\rangle^{1/2}\, . 
\end{eqnarray*}
Together with $G_{n,\beta}\t g\wstack{n\to\infty}{\lra}G_\beta \t g$ (cf. 
Step 1), this implies $G_{n,\beta}g_n\wstack{n\to\infty}{\lra}G_\beta g$. 
Similarly, we obtain 
\begin{eqnarray*}
&&\hspace{-.5cm}\left|\left(\langle G_{n,\beta}g_n\, , \, G_{n,\beta}g_n
\rangle_n-\langle G_\beta g\, , \, G_\beta g\rangle\right)-\left(\langle 
G_{n,\beta}\t g\, , \, G_{n,\beta}\t g\rangle_n-\langle G_\beta\t g\, , \, 
G_\beta\t g\rangle\right)\right| \\ 
&&\hspace{0.5cm}\le \left|\langle G_{n,\beta}(g_n+\t g)\, ,\, G_{n,\beta} 
(g_n-\t g)\rangle_n\right|+\left|\langle G_\beta (g+\t g)\, ,\, G_\beta (g 
-\t g)\rangle\right| \vphantom{\sum}\\ 
&&\hspace{0.5cm}\le {\textstyle\frac{1}{\beta^2}}\sup_{n\in {\Bbb N}}\langle 
g_n+\t g\, ,\, g_n+\t g\rangle_n^{1/2}\left(\langle g_n\, ,\, g_n\rangle_n-
2\langle g_n\, ,\, \t g\rangle_n+\langle \t g\, , \, \t g\rangle_n\right 
)^{1/2} \\ 
&&\hspace{1.0cm}+{\textstyle\frac{\ve}{\beta^2}}\langle g+\t g\, ,\, g+ 
\t g\rangle^{1/2}
\end{eqnarray*}
and thus 
\begin{eqnarray*}
&&\hspace{-.5cm}\limsup_{n\to\infty}\left|\left(\langle G_{n,\beta}g_n\, , 
\, G_{n,\beta}g_n\rangle_n-\langle G_\beta g\, , \, G_\beta g\rangle\right) 
-\left(\langle G_{n,\beta}\t g\, , \, G_{n,\beta}\t g\rangle_n-\langle G_\beta 
\t g\, , \, G_\beta\t g\rangle\right)\right| \\ 
&&\hspace{0.5cm}\le {\textstyle\frac{1}{\beta^2}} 
\sup_{n\in {\Bbb N}}\langle g_n+\t g\, , \, g_n+\t g\rangle_n^{1/2}\left 
\langle g-\t g\, , \, g-\t g\right\rangle^{1/2}+{\textstyle\frac{\ve} 
{\beta^2}}\langle g+\t g\, , \, g+\t g\rangle^{1/2} \\ 
&&\hspace{0.5cm}={\textstyle\frac{\ve}{\beta^2}}\left(\sup_{n\in {\Bbb N}} 
\langle g_n+\t g\, ,\, g_n+\t g\rangle_n^{1/2}+\left\langle g+\t g\, ,\, g 
+\t g\right\rangle^{1/2}\right)\, , 
\end{eqnarray*}
note that $\sup_{n\in {\Bbb N}}\langle g_n+\t g\, , \, g_n+\t g 
\rangle_n^{1/2}<\infty$ since $g_n+\t g\sstack{n\to\infty}{\lra}g+\t g$. 
Together with $G_{n,\beta}\t g\sstack{n\to\infty}{\lra}$ $G_\beta \t g$ 
(cf. Step 1), and $G_{n,\beta}g_n\wstack{n\to\infty}{\lra}G_\beta g$, 
this implies $G_{n,\beta}g_n\sstack{n\to\infty}{\lra}G_\beta g$. 
\medskip

\nid 
{\it Step 4 } Let $f\in L^2(E,\bnu)$, $f_n\in {\cal C}$, $n\in 
{\Bbb N}$, be a sequence $w$-converging to $f$, and let $\vp\in {\cal 
C}$. By the result of Step 1 and Proposition \ref{Proposition2.3} (c), 
we have 
\begin{eqnarray*}
\langle G_{n,\beta}'f_n\, ,\, \vp\rangle_n&=&\langle f_n\, ,\, G_{n, 
\beta}\vp\rangle_n \\ 
&\stack{n\to\infty}{\lra}&\langle f\, ,\, G_{\beta}\vp\rangle \\ 
&=&\langle G_{\beta}'f\, ,\, \vp\rangle \, . 
\end{eqnarray*}
This means nothing but $G_{n,\beta}'f_n\wstack{n\to\infty}{\lra}G_\beta' 
f$. Moreover, let $g\in L^2(E,\bnu)$ and $g_n\in {\cal C}$ $s$-converging 
to $g$ as $n\to\infty$. We have just proved that $G_{n,\beta}'g_n\wstack 
{n\to\infty}{\lra}G_\beta'g$ which now implies $G_{n,\beta}G_{n,\beta}'g_n 
\wstack{n\to\infty}{\lra}G_\beta G_\beta'g$. Again by Proposition 
\ref{Proposition2.3} (c), 
\begin{eqnarray*}
\langle G_{n,\beta}'g_n\, ,\, G_{n,\beta}'g_n\rangle_n&=&\langle g_n\, ,\, 
G_{n,\beta}G_{n,\beta}'g_n\rangle_n \\ 
&\stack{n\to\infty}{\lra}&\langle g\, ,\, G_{\beta}G_\beta'g\rangle \\ 
&=&\langle G_{\beta}'g\, ,\, G_{\beta}'g\rangle\, . 
\end{eqnarray*}
Thus, $G_{n,\beta}'g_n$ $s$-converges to $G_\beta'g$ as $n\to\infty$. 
\qed 
\begin{lemma}\label{Lemma2.8} 
Let $S_n$, $n\in {\Bbb N}$, be a sequence of bilinear forms pre-convergent 
to $S$. Condition (iii) of Lemma \ref{Lemma2.6} implies condition (iv). 
\end{lemma} 
Proof. Suppose we have condition (iii). Let us also assume that 
$\beta u_n-A_nu_n$ $w$-converges to $\beta u-Au$. We show that $u_n$ 
$w$-converges to $u$. 

For all $\psi\in {\cal C}$ and all sequences $\psi_n\in {\cal C}$, 
$n\in {\Bbb N}$, $s$-convergent to $\psi$, we have because of 
Proposition \ref{Proposition2.3} (c), 
\begin{eqnarray*} 
S_{n,\beta}(u_n,\psi_n)&=&\langle \beta u_n-A_nu_n\, , \, \psi_n 
\rangle_n \\ 
&\stack{n\to\infty}{\lra}&\langle \beta u-Au\, , \, \psi\rangle \\ 
&=&S_\beta(u,\psi)\, . 
\end{eqnarray*}
From condition (iii), it follows that 
\begin{eqnarray*} 
S_{n,\beta}(\psi_n,u_n)&\stack{n\to\infty}{\lra}&S_\beta(\psi,u)\, . 
\end{eqnarray*}
for all $\psi\in D(S)$ and all sequences $\psi_n\in D(S_n)\cap {\cal 
C}$, $n\in {\Bbb N}$, $s$-convergent to $\psi$ in the sense of 
condition (ii) in Definition \ref{Definition2.4}. In particular, we 
can take $g\in {\cal C}$, $\psi_n:=G_{n,\beta}g$, $n\in {\Bbb N}$, 
and $\psi:=G_{\beta}g$. For this, recall also condition (c3). 
According to Theorem \ref{Theorem2.7}, $\psi_n\sstack{n\to\infty} 
{\lra}\psi$. Furthermore, $\psi_n\in D(S_n)\cap {\cal C}$, $n\in 
{\Bbb N}$, and satisfies condition (ii) in Definition 
\ref{Definition2.4} since $\sup_{n\in {\Bbb N}}\langle A_n\psi_n\, , 
\, A_n\psi_n\rangle_n=\sup_{n\in {\Bbb N}}\langle\beta\psi_n-g\, ,\, 
\beta\psi_n-g\rangle_n<\infty$ and 
\begin{eqnarray*}
\limsup_{n\to\infty}S_n(\psi_n,\psi_n)=\lim_{n\to\infty}\langle 
g-\beta\psi_n\, ,\, \psi_n\rangle_n=\langle g-\beta \psi\, ,\, \psi 
\rangle=S(\psi ,\psi )\, .  
\end{eqnarray*}
Thus, 
\begin{eqnarray*}
\langle g,u_n\rangle_n=S_{n,\beta}(\psi_n,u_n)\stack{n\to\infty} 
{\lra}S_\beta(\psi,u)=\langle g,u\rangle
\end{eqnarray*}
which means that $u_n\wstack{n\to\infty}{\lra}u$. We have verified (iv). 
\qed
\medskip 

\nid
{\bf Remarks} (5) Lemmas \ref{Lemma2.6} and \ref{Lemma2.8} 
together show that under (i) and (ii) of Definition \ref{Definition2.4}, 
conditions (iii) and (iv) are equivalent. 
\medskip 

\nid
(6) By virtue of Remark (1), condition ($c3$), and Theorem \ref{Theorem2.7}, 
for $g\in {\cal C}$ and $g_n\in {\cal C}$ $s$-convergent to $g$ it holds that 
\begin{eqnarray}\label{2.7}
\langle G_{n,\beta}g_n - G_\beta g\, ,\, G_{n,\beta}g_n - G_\beta g\rangle_n 
\stack{n\to\infty}{\lra} 0\, . 
\end{eqnarray}
For $f\in L^2(E,\bnu)$ and $f_n\in {\cal C}$ $w$-convergent to $f$, from 
Theorem \ref{Theorem2.7} and Proposition \ref{Proposition2.3} (e) it follows 
that 
\begin{eqnarray*} 
\liminf_{n\to\infty}\langle G_{n,\beta}f_n\, ,\, G_{n,\beta}f_n\rangle_n\ge 
\langle G_\beta f\, ,\, G_\beta f\rangle 
\end{eqnarray*} 
as well as 
\begin{eqnarray*} 
\liminf_{n\to\infty}\langle G_{n,\beta}'f_n\, ,\, G_{n,\beta}'f_n\rangle_n 
\ge\langle G_\beta'f\, ,\, G_\beta'f\rangle\, .  
\end{eqnarray*} 
In order to prove convergence of forms $S_n$, $n\in {\Bbb N}$, to a 
form $S$ let us introduce one more condition: 
\begin{itemize}
\item[(c4)] $D(S)\subseteq {\cal C}$ in the sense that for every 
$\vp\in D(S)$, there is a $u\in {\cal C}$ with $\psi=u$ ${\bnu}$-a.e. 
Moreover, in this sense $\{A\vp:\vp\in D(S)\}\subseteq {\cal C}$. 
\end{itemize} 

In the next proposition we formulate conditions under which (i)-(iii) 
of Definition \ref{Definition2.4} become necessary. For this recall 
Remark (3).  
\begin{proposition}\label{Proposition2.9} 
Suppose that conditions (c1) -- (c4) are satisfied. Let $(G_{\beta})_{ 
\beta\ge 0}$ be the resolvent of a strongly continuous contraction 
semigroup $(T_{t})_{t\ge 0}$ on $L^2(E,\bnu)$ and let $(G_{n,\beta})_{ 
\beta\ge 0}$ be the resolvent of a strongly continuous contraction 
semigroup $(T_{n,t})_{t\ge 0}$ on $L^2(E,\bnu_n)$, $n\in {\Bbb N}$. 
Suppose we have the following. 
\begin{itemize} 
\item[(i)] $G_{n,\beta}g_n$ $s$-converges to $G_\beta g$ as $n\to 
\infty$ for every $g\in {\cal C}$, every sequence $g_n\in D(S_n)\cap 
{\cal C}$ with $g_n\sstack{n\to\infty}{\lra}g$, and $\beta >0$. 
Furthermore, $G_{n,\beta}'g$ $s$-converges to $G_\beta 'g$ as $n\to 
\infty$ for every $g\in {\cal C}$ and $\beta >0$.
\item[(ii)] $\psi\in D(S_n)\cap {\cal C}=D(A_n)\cap {\cal C}$ implies 
$\psi\in D(A'_n)$. 
\item[(iii)] Let $D(S_{n_k})\cap {\cal C}\ni\vp_{n_k}\wstack{n\to 
\infty}{\lra}\vp\in L^2(E,\bnu)$ such that $\sup_{k\in {\Bbb N}}\left 
\langle A_{n_k}\vp_{n_k}\, ,\, A_{n_k}\vp_{n_k}\right\rangle_{n_k}< 
\infty$. Then $\vp\in D(S)$ and $\sup_{k\in {\Bbb N}}\left\langle 
A'_{n_k}\vp_{n_k}\, ,\, A'_{n_k}\vp_{n_k}\right\rangle_{n_k}<\infty$.
\end{itemize} 
Then the forms $S_n$, $n\in {\Bbb N}$, associated with $(T_{n,t})_{ 
t\ge 0}$ converge to the form $S$ associated with $(T_{t})_{t\ge 0}$ 
as $n\to\infty$ in the sense of Definition \ref{Definition2.4}. 
\end{proposition} 
Proof. {\it Step 1 } In the first two steps, we verify condition (i) 
of Definition \ref{Definition2.4}. For a clearer presentation, we 
ignore subsequences.  Let us introduce bilinear forms $S^{(\beta)} 
(w,w):=\beta\left\langle w-\beta G_{\beta}w\, , \, w\right\rangle$, 
$w\in D(S)$, and $S^{(\beta)}_n(w_n,w_n):=\beta\left\langle w_n-\beta 
G_{n,\beta}w_n\, , \, w_n\right\rangle_n$, $w_n\in D(S_n)\cap {\cal 
C}$, $n\in {\Bbb N}$, which can be considered a counterpart to the 
Deny-Yosida approximation in Dirichlet form theory. Classical 
semigroup theory says that $\beta\left(w-\beta G_{\beta}w\right)\stack 
{\beta\to\infty}{\lra}-Aw$ in $L^2(E,\bnu)$ and therefore 
\begin{eqnarray*} 
S^{(\beta)}(w,w)=\beta\langle w-\beta G_{\beta}w\, ,\, w\rangle\stack 
{\beta\to\infty}{\lra}S(w,w)\vphantom{\sum}\, ,\quad w\in D(S). 
\end{eqnarray*}
and similarly $\lim_{\beta\to\infty}S_n^{(\beta)}(w_n,w_n)=S_n(w_n, 
w_n)$, $w_n\in D(S_n)\cap {\cal C}$, $n\in {\Bbb N}$. Differentiating 
with respect to $\beta$ and applying the resolvent identity we get 
\begin{eqnarray*} 
\left|\frac{d}{d\beta}S_n^{(\beta)}(w_n,w_n)\right|&=&\left|\frac{d} 
{d\beta}\beta\langle w_n-\beta G_{\beta}w_n\, ,\, w_n\rangle_n\right| 
 \\ 
&=&\left|\left\langle w_n-2\beta G_{n,\beta}w_n+\beta^2G_{n,\beta}^2 
w_n\, ,\, w_n\right\rangle_n\vphantom{\dot{f}}\right|\vphantom{\frac 
{d}{d\beta}}\vphantom{\frac{d}{d\beta}} \\ 
&=&\left|\left\langle (A_nG_{n,\beta})^2w_n\, ,\, w_n\right\rangle_n 
\vphantom{\dot{f}}\right|\vphantom{\frac{d}{d\beta}} \\ 
&=&\left|\left\langle G_{n,\beta}A_nG_{n,\beta}w_n\, ,\, A'_nw_n\right 
\rangle_n\vphantom{\dot{f}}\right|\vphantom{\frac{d}{d\beta}} \\ 
&=&\left|\left\langle G_{n,\beta}^2A_nw_n\, ,\, A'_nw_n\right\rangle_n
\vphantom{\dot{f}}\right|\vphantom{\frac{d}{d\beta}} \\ 
&\le&\frac{1}{\beta^2}\left\langle A_nw_n\, ,\, A_nw_n\right 
\rangle_n^{1/2}\langle A_n'w_n\, ,\, A_n'w_n\rangle_n^{1/2}\, ,\quad 
w_n\in D(S_n)\cap {\cal C}, \vphantom{\frac{d}{d\beta}}
\end{eqnarray*}
where, for the last three lines we have applied hypothesis (ii). It 
follows now from hypothesis (iii) that $\lim_{\beta\to\infty}S_n^{ 
(\beta)}(w_n,w_n)=S_n(w_n,w_n)$ uniformly in $n\in {\Bbb N}$ whenever 
$\sup_{n\in {\Bbb N}}\left\langle A_n\vp_n\, ,\right.$ $\left.A_n 
\vp_n\right\rangle_n<\infty$ and $D(S_n)\cap {\cal C}\ni w_n\wstack{n 
\to\infty}{\lra}w\in L^2(E,\bnu)$. 
\medskip

\nid 
{\it Step 2 }  Let us complete the verification of condition (i) of 
Definition \ref{Definition2.4}. Let $D(S_n)\cap {\cal C}\ni\vp_n 
\wstack{n\to\infty}{\lra}\vp\in {\cal C}$ and recall (c4). We can 
decompose 
\begin{eqnarray*} 
&&\hspace{-.5cm}S_n^{(\beta)}(\vp_n,\vp_n)-S^{(\beta)}(\vp,\vp)= 
\beta\left\langle\vp_n-\beta G_{n,\beta}\vp_n\, ,\, \vp_n\right 
\rangle_n-\beta\left\langle\vp-\beta G_{\beta}\vp\, ,\, \vp 
\right\rangle\vphantom{sum} \\ 
&&\hspace{.5cm}=\beta\langle\vp_n-\beta G_{n,\beta}\vp\, ,\, \vp 
\rangle_n+\beta\langle\vp-\beta G_{\beta}\vp\, ,\, \vp_n\rangle_n 
\vphantom{\sum} \\ 
&&\hspace{1cm}-\beta\langle\vp-\beta G_{\beta}\vp\, ,\, \vp 
\rangle_n-\beta\langle\vp-\beta G_{\beta}\vp\, ,\, \vp\rangle 
\vphantom{\sum} \\ 
&&\hspace{1cm}+\beta^2\langle G_{\beta}\vp-G_{n,\beta}\vp-G_{n, 
\beta}'\vp,\, \vp_n-\vp\rangle_n\vphantom{\sum} \\ 
&&\hspace{1cm}+\beta\langle\vp_n-\vp-\beta G_{n,\beta}(\vp_n-\vp) 
\, ,\, \vp_n-\vp\rangle_n\, . \vphantom{\sum} 
\end{eqnarray*}
Let us analyze the items on the right-hand side. 
\begin{itemize}
\item[(1)] $\langle\vp_n-\beta G_{n,\beta}\vp\, , \, \vp\rangle_n 
\stack{n\to\infty}{\lra}\langle \vp-\beta G_{\beta}\vp\, , \, \vp 
\rangle$ because of $\vp_n-\beta G_{n,\beta}\vp\wstack{n\to\infty}
{\lra}\vp-\beta G_{\beta}\vp$, cf. (i) of this proposition, 
\item[(2)] $\langle\vp-\beta G_{\beta}\vp\, , \, \vp_n\rangle_n 
\stack{n\to\infty}{\lra}\langle \vp-\beta G_{\beta}\vp\, , \, \vp 
\rangle$ because of $\vp-\beta G_{\beta}\vp\in {\cal C}$ (cf. 
condition (c3) and Lemma 2.1 (b)) and $\vp_n\wstack{n\to\infty} 
{\lra}\vp$,  
\item[(3)] $\langle\vp-\beta G_{\beta}\vp\, , \, \vp\rangle_n 
\stack{n\to\infty}{\lra}\langle \vp-\beta G_{\beta}\vp\, , \, \vp 
\rangle$ because of condition (c3), Lemma \ref{Lemma2.1} (b) and, 
thus, $\vp-\beta G_{\beta}\vp\in {\cal C}$, and Lemma \ref{Lemma2.1} 
(d), 
\item[(4)] $\langle G_{\beta}\vp-G_{n,\beta}\vp-G_{n,\beta}'\vp , \, 
\vp_n -\vp \rangle_n\stack{n\to\infty}{\lra}0$ because of $G_{\beta} 
\vp-G_{n,\beta}\vp-G_{n,\beta}'\vp\sstack{n\to\infty}{\lra}-G_\beta' 
\vp$ (cf. (i) of this proposition), $\vp_n -\vp\wstack{n\to\infty} 
{\lra}0$, and Proposition \ref{Proposition2.3} (c), 
\item[(5)] $\langle \vp_n-\vp-\beta G_{n,\beta}(\vp_n-\vp)\, , \, 
\vp_n-\vp\rangle_n\ge 0$ because of Lemma 2.1 (b) in \cite{Lo13}. 
\end{itemize}
Therefore, $\liminf_{n\to\infty}S_n^{(\beta)}(\vp_n,\vp_n)\ge 
S^{(\beta)}(\vp,\vp)$, $\beta>0$. With the result of Step 1, we get 
$\liminf_{n\to\infty}$ $S_n(\vp_n,\vp_n)\ge S(\vp,\vp)$ which, by 
condition (iii) of the present proposition and (c4), is condition (i) 
of Definition \ref{Definition2.4}. 
\medskip

\nid 
{\it Step 3 } Let us verify condition (ii) of Definition 
\ref{Definition2.4}. Among others things, we will use an idea from 
the proof of \cite{Mo94}, Theorem 2.4.1, part (jj). Let $\psi\in 
D(S)$ and recall that, because of (c4), we have $\psi\in {\cal C}$. 
Moreover let $\Psi:=\psi-A\psi$, i. e. $\Psi\in {\cal C}$ by (c4) and 
Lemma \ref{Lemma2.1} (b). Furthermore, $\psi=G_1\Psi$. Let $\t\psi_n: 
=G_{n,1}\Psi$, i. e. $\t\psi_n\in {\cal C}$, $n\in {\Bbb N}$, by (c3). 
By condition (i) of the present proposition we have $\t\psi_n=G_{n,1} 
\Psi\sstack{n\to\infty}{\lra}G_1\Psi=\psi$ and therefore $G_{n,\beta} 
\t\psi\sstack{n\to\infty}{\lra}G_\beta\psi$. According to Proposition 
\ref{Proposition2.3} (c) this means 
\begin{eqnarray*}
S_n^{(\beta)}(\t\psi_n,\t\psi_n)&&\hspace{-.5cm}=\beta\langle\t\psi_n 
-\beta G_{n,\beta}\t\psi_n\, ,\, \t\psi_n\rangle_n \\ 
&&\hspace{-1cm}\stack{n\to\infty}{\lra}\beta\langle\psi-\beta G_\beta 
\psi\, ,\, \psi\rangle=S^{(\beta)}(\psi,\psi)\, ,\quad\beta >0. 
\end{eqnarray*}

On the other hand, it holds that $S^{(\beta )}(\psi,\psi)\stack{\beta 
\to\infty}{\lra}S(\psi,\psi)$, cf. Step 1. Thus, there exists a sequence 
$\beta_n$, $n\in {\Bbb N}$, with $\beta_n\stack{n\to\infty}{\lra}\infty$ 
such that 
\begin{eqnarray}\label{2.8} 
\lim_{n\to\infty}S_n^{(\beta_n)}(\t\psi_n,\t\psi_n)=S(\psi,\psi)\, .  
\end{eqnarray}
Because of 
\begin{eqnarray}\label{2.9} 
S_n^{(\beta_n)}(\t\psi_n,\t\psi_n)&&\hspace{-.5cm}=S_n(\beta_n G_{n, 
\beta_n}\t\psi_n\, ,\, \beta_n G_{n,\beta_n}\t\psi_n)\nonumber \\ 
&&\hspace{-.0cm}+\beta_n\langle\t\psi_n-\beta_n G_{n,\beta_n}\t\psi_n 
\, ,\, \t\psi_n-\beta_nG_{n,\beta_n}\t\psi_n\rangle_n 
\end{eqnarray}
and $S_n(\beta_n G_{n,\beta_n}\t\psi_n\, ,\, \beta_n G_{n,\beta_n}\t 
\psi_n)\ge 0$ it follows from (\ref{2.8}) that 
\begin{eqnarray}\label{2.10}  
\langle\t\psi_n-\beta_nG_{n,\beta_n}\t\psi_n\, ,\, \t\psi_n-\beta_n 
G_{n,\beta_n}\t\psi_n\rangle_n\stack{n\to\infty}{\lra}0  
\end{eqnarray}
which implies 
\begin{eqnarray*} 
|\langle\t\psi_n-\beta_nG_{n,\beta_n}\t\psi_n\, ,\, \vp\rangle_n|&& 
\hspace{-.5cm}\le\langle\t\psi_n-\beta_nG_{n,\beta_n}\t\psi_n\, ,\, 
\t\psi_n-\beta_n G_{n,\beta_n}\t\psi_n\rangle_n^{1/2}\langle\vp\, ,\, 
\vp\rangle_n^{1/2} \\ 
&&\hspace{-1cm}\stack{n\to\infty}{\lra}0\, ,\quad\vp\in {\cal C}. 
\end{eqnarray*}
We recall that ${\cal C}$ is linear by Lemma \ref{Lemma2.1} (b). 
Thus $\beta_nG_{n,\beta_n}\t\psi_n\wstack{n\to\infty}{\lra}\psi$ and 
(\ref{2.10}) provides  
\begin{eqnarray*} 
\beta_nG_{n,\beta_n}\t\psi_n\sstack{n\to\infty}{\lra}\psi\, .  
\end{eqnarray*}
Relations (\ref{2.8}) and (\ref{2.9}) imply now 
\begin{eqnarray*} 
\limsup_{n\to\infty}S_n(\beta_nG_{n,\beta_n}\t\psi_n,\beta_n G_{n, 
\beta_n}\t\psi_n)\le S(\psi,\psi)\, . 
\end{eqnarray*}
In addition, by $\psi_n\in D(A_n)$, $n\in {\Bbb N}$, and $\Psi\in {\cal C} 
\subseteq{\cal D}$, 
\begin{eqnarray*} 
&&\hspace{-.5cm}\sup_{n\in {\Bbb N}}\left\langle A_n\beta_nG_{n,\beta_n} 
\t\psi_n\, ,\, A_n\beta_nG_{n,\beta_n}\t\psi_n\right\rangle_n=\sup_{n\in 
{\Bbb N}}\left\langle\beta_nG_{n,\beta_n}A_n\t\psi_n\, ,\, \beta_nG_{n, 
\beta_n}A_n\t\psi_n\right\rangle_n \\ 
&&\hspace{.5cm}\le\sup_{n\in {\Bbb N}}\left\langle A_n\t\psi_n\, ,\, A_n 
\t\psi_n\right\rangle_n=\sup_{n\in {\Bbb N}}\left\langle\t\psi_n-\Psi\, ,\, 
\t\psi_n-\Psi\right\rangle_n\le 4\sup_{n\in {\Bbb N}}\langle\Psi,\Psi 
\rangle_n<\infty\, . 
\end{eqnarray*}

Choosing $\psi_n:=\beta_n G_{n,\beta_n}\t\psi_n$, $n\in {\Bbb N}$, this 
verifies condition (ii) of Definition \ref{Definition2.4}. 
\medskip

\nid 
{\it Step 4 } In order to verify condition (iii) of Definition 
\ref{Definition2.4}, let $\beta >0$, $\psi\in {\cal C}$, and let $u_n 
\in D(S_n)$, $n\in {\Bbb N}$, be a $w$-convergent sequence with the 
properties mentioned in Lemma \ref{Lemma2.6} and $u\in D(S)$. We assume 
that $\lim_{n\to\infty}S_{n,\beta}(u_n,\psi_n)=S_\beta (u,\psi)$ for 
all sequences $\psi_n\in {\cal C}$ $s$-converging to $\psi$. Recalling 
hypothesis (i) it follows that for $g\in {\cal C}$, $\hat{\psi}_n:=G_{ 
n,\beta}'g$, $n\in {\Bbb N}$, and $\hat{\psi}:=G_\beta 'g$ we have 
\begin{eqnarray*} 
\lim_{n\to\infty}\langle u_n\, , \, g \rangle_n&=&\lim_{n\to\infty}S_{n, 
\beta}(u_n,\hat{\psi}_n) \\ 
&=&S_\beta (u,\hat{\psi}) \\ 
&=&\langle u\, , \, g \rangle\,.
\end{eqnarray*}
This means that $u_n$ $w$-converges to $u$. For this, also note that 
because of condition (c3), $\hat{\psi},\hat{\psi}_n\in {\cal C}$, $n\in 
{\Bbb N}$. Let now $\psi\in D(S)$ and $\psi_n\in D(S_n)\cap {\cal C}$, $n 
\in {\Bbb N}$, be a sequence $s$-convergent to $\psi$ in the sense of 
Definition \ref{Definition2.4} (ii). Assume $\lim_{n\to\infty}S_{n,\beta} 
(u_n,\psi_n)=S_\beta (u,\psi)$. Lemma \ref{Lemma2.5} yields $\lim_{n\to 
\infty}S_{n,\beta}(\psi_n,u_n)=S_\beta (\psi ,u)$. We have verified (iii) 
of Definition \ref{Definition2.4}. 
\qed 
\medskip 

\nid 
We conclude this subsection with the proof of $s$-convergence of the 
associated semigroups. For this, let us introduce the following condition. 
\begin{itemize} 
\item[(c5)] ${\cal T}:=\{T_t g:g\in {\cal C},\ t >0\}\subseteq {\cal 
C}$, ${\cal T}_n:=\{T_{n,t}g:g\in {\cal C},\ t >0\}\subseteq {\cal C}$, 
${\cal T}':=\{T'_t g:g\in {\cal C},\ t >0\}\subseteq {\cal C}$, and  
${\cal T}'_n:=\{T'_{n,t}g:g\in {\cal C},\ t >0\}\subseteq {\cal C}$, 
$n\in {\Bbb N}$, in the sense of condition (c3).  
\end{itemize} 
\begin{theorem}\label{Theorem2.10} 
Suppose that (c1),(c2),(c3),(c5) are satisfied. Then, for all $g\in {\cal 
C}$ and \linebreak $\beta >0$, $G_{n,\beta}g\sstack{n\to\infty}{\lra}G_\beta g$ iff 
$\ T_{n,t}g\sstack{n\to\infty}{\lra} T_tg$ for all $t>0$ and $G_{n,\beta}' 
g\sstack{n\to \infty}{\lra} G_\beta 'g$ iff \linebreak $T'_{n,t}g\sstack{n\to\infty} 
{\lra} T'_tg$ for all $t>0$. 
\end{theorem} 
Proof. In Steps 1 and 2, we demonstrate that $G_{n,\beta}g\sstack{n\to 
\infty}{\lra}G_\beta g$ implies $T_{n,t}g\sstack{n\to\infty}{\lra} T_tg$ 
and that and $G'_{n,\beta}g\sstack{n\to\infty}{\lra}G'_\beta g$ implies 
$T'_{n,t}g\sstack{n\to\infty}{\lra}T'_tg$. In Step 3 we verify the 
converse. 
\medskip 

\nid 
{\it Step 1 } In this step, let us show that $T_{n,t}f\wstack{n 
\to\infty}{\lra}T_tf$ for all $f\in {\cal C}$. 
For well-definiteness, recall condition (c5). Since both, ${\cal C}$ and 
$D(A^2)$, are dense in $L^2(E,{\bnu})$ it is sufficient to verify this 
claim for $f=(G_\beta)^2 h$, $h\in {\cal C}$. Set $g=G_\beta h$. With 
\cite{P83}, Lemma 4.1 of Chapter 3, we have 
\begin{eqnarray*} 
&&\hspace{-.5cm}\langle T_{n,t}f-T_tf\, ,\, \vp\rangle_n=\langle T_{n,t} 
G_\beta g-T_tG_\beta g\, ,\, \vp\rangle_n \\ 
&&\hspace{.5cm}=\langle T_{n,t}G_\beta g-T_{n,t}G_{n,\beta}g\, ,\, \vp 
\rangle_n+\langle G_{n,\beta}T_{n,t}g-G_{n,\beta}T_tg\, ,\, \vp\rangle_n 
+\langle G_{n,\beta}T_tg-G_{\beta}T_tg\, ,\, \vp\rangle_n \\ 
&&\hspace{.5cm}=\langle T_{n,t}G_\beta g-T_{n,t}G_{n,\beta}g\, ,\, \vp 
\rangle_n-{\textstyle\left\langle\int_0^t T_{n,t-s}(G_\beta -G_{n,\beta}) 
T_{s}h\, ds\, ,\, \vp \right\rangle_n} \\ 
&&\hspace{1cm}+\langle G_{n,\beta}T_tg-G_{\beta} T_tg\, ,\, \vp\rangle_n 
\, .
\end{eqnarray*}
It follows from the Schwarz inequality, (\ref{2.7}), and contractivity 
of $T_{n,t}$ in $L^2(E,\bnu_n)$ that the first item of the right-hand side 
tends to zero. That the third item of the right-hand side tends to zero is 
a consequence of the Schwarz inequality and (\ref{2.7}). Thus, it remains 
to demonstrate that $\langle\int_0^t T_{n,t-s}(G_\beta -G_{n,\beta}) T_{s} 
h\, ds\, ,\, \vp \rangle_n$ $\stack{n\to\infty}{\lra}0$. But using again 
Schwarz' inequality, and contractivity of all the $T_{n,t-s}$, this follows 
from 
\begin{eqnarray}\label{2.11}
{\textstyle\left|\left\langle\int_0^t T_{n,t-s}(G_\beta -G_{n,\beta})T_{s}h 
\, ds\, , \, \vp \right\rangle_n\right|}\nonumber \\ 
&&\hspace{-5cm}\le \ {\textstyle\int_0^t\left\langle (G_\beta -G_{n, 
\beta})T_{s}h\, , \, (G_\beta -G_{n,\beta})T_{s}h\right\rangle_n^{1/2} 
\left\langle \vp \, , \, \vp \right\rangle_n^{1/2}\, ds}\, , 
\end{eqnarray}
condition (c5), relation (\ref{2.7}), and dominated convergence. 
\medskip 

\nid 
{\it Step 2 } Now, let us prove that $T_{n,t}f\sstack{n\to\infty}{\lra}T_tf$ 
for all $f\in {\cal C}$. Again, we choose $f=(G_\beta)^2h$ where $h\in {\cal 
C}$, set $g=G_\beta h$, and decompose 
\begin{eqnarray*} 
&&\hspace{-.5cm}T_{n,t}f-T_tf \\ 
&&\hspace{.5cm}=\left(T_{n,t}G_\beta g-T_{n,t}G_{n,\beta} g\right)-\left({ 
\textstyle\int_0^t T_{n,t-s}(G_\beta -G_{n,\beta})T_{s}h\, ds}\right)+\left( 
G_{n,\beta}T_tg-G_\beta T_tg\right)\, . 
\end{eqnarray*}
That the $L^2(E,{\bnu}_n)$-norms of the first and the third item tend to 
zero as $n\to\infty$ follows from the contractivity of the semigroups and 
relation (\ref{2.7}). The arguments used already in Step 1, (\ref{2.11}), 
lead to 
\begin{eqnarray*} 
\left\langle{\textstyle\int_0^t T_{n,t-s}(G_\beta -G_{n,\beta})T_{s}h\, 
ds}\, , \, {\textstyle\int_0^t T_{n,t-s}(G_\beta -G_{n,\beta})T_{s}h\, ds} 
\right\rangle_n\stack{n\to\infty}{\lra}0\, . 
\end{eqnarray*}
For the dual operators, we recall that also $(T_{n,t}')_{t\ge 0}$ is a 
strongly continuous contraction semigroup in $L^2(E,\bnu_n)$, $n\in {\Bbb 
N}$, and $(T_t')_{t\ge 0}$ is a strongly continuous contraction semigroup 
in $L^2 (E,\bnu)$, cf. \cite{P83} Subsection 1.10, especially Corollary 
10.6. This means we can conclude that $T_{n,t}'f\sstack{n\to\infty}{\lra} 
T_t'f$ for all $f\in {\cal C}$ the same way as we did it to show $T_{n,t} 
f\sstack{n\to\infty}{\lra}T_tf$. 
\medskip 

\nid
{\it Step 3 } This follows from $G_\beta f=\int_0^\infty e^{-\beta t} 
T_tf\, dt$, $G_{n,\beta} f=\int_0^\infty e^{-\beta t}T_{n,t}f\, dt$, $n 
\in {\Bbb N}$, $f\in {\cal C}$, the same relations for the dual operators, 
and dominated convergence. \qed 
\medskip 

\nid
{\bf Remark} (7) From Theorems \ref{Theorem2.7} and \ref{Theorem2.10} it 
follows that for $g\in {\cal C}$ it holds that $\langle T_{n,t}g-T_t g\, 
,\, T_{n,t}g-T_t g\rangle_n\stack{n\to\infty}{\lra} 0$ and $\langle T_{n, 
t}'g-T_t' g\, ,\, T_{n,t}'g-T_t' g\rangle_n\stack{n\to\infty}{\lra} 0$.  
Even more important is the following observation. By the contractivity 
of the semigroups $(T_{n,t})_{t\ge 0}$ we have for ${\cal C}\ni g_n 
\sstack{n\to\infty}{\lra}g\in {\cal C}$ the limits $\langle T_{n,t}g_n 
-T_{n,t}g\, ,\, T_{n,t}g_n-T_{n,t}g\rangle_n\stack{n\to\infty}{\lra}0$ 
and $\langle G_{n,\beta}g_n-G_{n,\beta}g\, ,\, G_{n,\beta}g_n-G_{n,\beta} 
g\rangle_n\stack{n\to\infty}{\lra}0$. 

Under the hypotheses of Theorem \ref{Theorem2.10} this says $G_{n,\beta} 
g_n\sstack{n\to\infty}{\lra}G_\beta g$ iff $G_{n,\beta}g\sstack{n\to 
\infty}{\lra}G_\beta g$ iff $\ T_{n,t}g\sstack{n\to\infty}{\lra} T_tg$ iff 
$\ T_{n,t}g_n\sstack{n\to\infty}{\lra} T_tg$ for all $\beta>0$ and $t>0$. 
The same holds for the dual operators.  

\subsection{Mosco Type Convergence of Non-Positive Non-Symmetric \\ Forms}

For the remainder of this section, let us drop the assumption that the 
semigroups $(T_{n,t})_{t\ge 0}$, $n\in {\Bbb N}$, and $(T_t)_{t\ge 0}$ 
are contractive. Anything else for the semigroups remains as introduced 
in Section 1. As a consequence, we cannot state positivity of the 
associated bilinear forms. We are interested in substitutes for Theorem 
\ref{Theorem2.7}. For this, we collect almost everything still necessary 
for the remainder of the section in the following condition.  
\begin{itemize}
\item[(c6)] 
\begin{itemize} 
\item[$(i)$] $\1\in D(A')$ and $\1\in D(A'_n)$, $n\in {\Bbb N}$, and $\sup_{ 
n\in {\Bbb N}}\|A_n'\1\|_{L^\infty (E,\sbnu_n)}<\infty$. 
\item[$(ii)$] $T_t'\1\in L^\infty (E,\bnu)$ and the limit $A'\1=\lim_{t\to 
0}\frac1t(T_t'\1-\1)$ exists in $L^\infty (E,\bnu)$. 
\item[$(iii)$] There exist $N_n\in {\cal B}(E_n)$ with $\bnu_n(N_n)\stack{n 
\to\infty}{\lra}0$ such that for $\vp\in {\cal C}$ there is $\Phi_n\equiv 
\Phi_n(\vp)\in {\cal C}$ with $\Phi_n=A_n'\1\cdot\vp$ on $E_n\setminus N_n$. 
Furthermore, $\Phi_n\sstack{n\to\infty}{\lra}A'\1\cdot\vp$. Alternatively 
assume $A'_n\1 \in {\cal V}$, $n\in {\Bbb N}$, and $A'_n\1\cdot\vp\sstack{n 
\to\infty}{\lra}A'\1\cdot\vp$ for all $\vp\in {\cal C}$. 
\item[$(iv)$] For $n\in {\Bbb N}$, there exists a set $D_n\subseteq D(S_n) 
\cap L^\infty (E,\bnu_n)$ which is dense in $D(S_n)$ with respect to the 
norm $\|f\|_{D_n}:=\left(\langle f\, , \, f\rangle_n+\langle A_nf\, , \, 
A_nf\rangle_n\right)^{1/2}$. 
\end{itemize}
\end{itemize}
Define $D(\hat{S}_n):=D(S_n)$, $n\in {\Bbb N}$, $D(\hat{S}):=D(S)$, and 
\begin{eqnarray}\label{2.12}
\hat{S}_n(u_n,v_n):=S_n(u_n,v_n)+{\textstyle\frac12}\langle A_n'\1\cdot u_n 
\, , \, v_n\rangle_n\, , \quad u_n\in D(\hat{S}_n),\ v_n\in L^2(E,\bnu_n),\ 
n\in {\Bbb N}, \ \ 
\end{eqnarray} 
and 
\begin{eqnarray}\label{2.13}
\hat{S}(u,v):=S(u,v)+{\textstyle\frac12}\langle A'\1\cdot u\, , \, v\rangle 
\, , \quad u\in D(\hat{S}),\ v\in L^2(E,\bnu). 
\end{eqnarray} 
\begin{lemma}\label{Lemma2.11} 
Suppose (c6). We have $\hat{S}_n(u_n,u_n)\ge 0$, $u_n\in D(S_n)$, $n\in 
{\Bbb N}$, and $\hat{S}(u,u)\ge 0$, $u\in D(S)$. 
\end{lemma} 
Proof. {\it Step 1 } We prove the claim for $\hat{S}$. Note that $T_t'\1\in 
L^\infty(E,\bnu)$ (cf. (c6$(ii)$)) and observe 
\begin{eqnarray*}
\|T_tu\|_{L^1(E,\sbnu)}\le \|T_t|u|\|_{L^1(E,\sbnu)}=\langle T_t'\1\, , \, 
|u|\rangle\, ,  \quad u\in L^2(E,\bnu).   
\end{eqnarray*} 
This implies 
\begin{eqnarray*}
\langle T_tu\, , \, T_tu\rangle\le\|T_tu^2\|_{L^1(E,\sbnu)}\le\int T_t'\1 
\cdot u^2\, d\bnu\, , \quad u\in L^2(E,\bnu),\ t\ge 0,  
\end{eqnarray*} 
the first inequality in this line by the association with a transition 
probability function. Consequently, 
\begin{eqnarray*}
\langle T_tu\, , \, u\rangle\le&&\hspace{-.5cm}\left(\int T_t'\1\cdot u^2\, 
d\bnu\right)^{1/2}\cdot\langle u\, , \, u\rangle^{1/2} \\ 
=&&\hspace{-.5cm}\langle u\, , \, u\rangle+\langle u\, , \, u\rangle^{1/2} 
\left(\left(\int T_t'\1\cdot u^2\, d\bnu\right)^{1/2}-\langle u\, , \, u 
\rangle^{1/2}\right) 
\end{eqnarray*} 
and therefore 
\begin{eqnarray*}
\ -\frac{\langle u\, , \, u\rangle^{1/2}}{\langle u\, , \, u\rangle^{1/2} 
+\left(\int T_t'\1\cdot u^2\, d\bnu\right)^{1/2}}\, {\textstyle\left\langle 
\frac1t\left(T_t'\1-\1\right)u\, , \, u\right\rangle}\le {\textstyle\frac1t} 
\langle u-T_tu\, , \, u\rangle\, , \quad u \in L^2(E,\bnu). 
\end{eqnarray*} 
Recalling $\1\in D(A')$ (cf. (c6$(i)$)) and (c6$(ii)$), and letting $t\to 
0$, it turns out that 
\begin{eqnarray*}
-{\textstyle\frac12}\langle A'\1\cdot u\, , \,u\rangle\le S(u,u)\, , \quad 
u\in D(S). 
\end{eqnarray*} 
{\it Step 2 } We prove the claim for $\hat{S}_n$. For this, we choose 
$u_n\in D_n$ (cf. (c6$(iv)$)) and proceed as in Step 1. We arrive at 
\begin{eqnarray*}
\ -\frac{\langle u_n\, , \, u_n\rangle_n^{1/2}}{\langle u_n\, , \, u_n 
\rangle_n^{1/2}+\left(\int T_{n,t}'\1\cdot u_n^2\, d\bnu_n\right)^{1/2}}\, 
{\textstyle\left\langle\frac1t\left(T_{n,t}'\1-\1\right)u_n\, , \, u_n 
\right\rangle_n}\le {\textstyle\frac1t}\langle u_n-T_{n,t}u_n\, , \, u_n 
\rangle_n\, . 
\end{eqnarray*} 
Letting again $t\to 0$, we obtain $-{\textstyle\frac12}\langle A_n'\1\cdot 
u_n\, , \,u_n\rangle_n\le S_n(u_n,u_n)$ for all $u_n\in D_n$. Let us finally 
mention that, in contrast to Step 1, we do not require $T'_{n,t}\1\in 
L^\infty (E,\bnu_n)$. We compensate this by requiring $u_n\in D_n\subseteq 
L^\infty (E,\bnu_n)$. But with (c6$(iv)$), we get finally $-{\textstyle 
\frac12}\langle A_n'\1\cdot u_n\, , \,u_n\rangle_n\le S_n(u_n,u_n)$ for all 
$u_n\in D(S_n)$. 
\qed 
\medskip

\nid
Let us assume that there are Markov processes associated with the 
semigroups $(T_{n,t})_{t\ge 0}$, $n\in {\Bbb N}$, and $(T_t)_{t\ge 
0}$: For $n\in {\Bbb N}$, let $X^n=((X^n_t)_{t\ge 0},(P^n_\mu)_{\mu 
\in E_n})$ be a process taking values in $E_n$ which corresponds to 
the semigroup $(T_{n,t})_{t\ge 0}$ and the form $S_n$. Here, $E_n$ is 
the subset of $E$ specified in Subsection 2.1. Furthermore, 
let $X=((X_t)_{t\ge 0},(P_\mu)_{\mu\in E})$ be a process associated 
with the semigroup $(T_t)_{t\ge 0}$ and the form $S$ which takes 
values in some subset of $E$. Suppose that the paths of the processes 
$X^n$, $n\in {\Bbb N}$, and $X$ are cadlag. For $\beta >0$, introduce 
$G_{n,\beta}g_n:=\int_0^\infty e^{-\beta t}T_{n,t}g_n\, dt$, $g_n\in 
L^\infty (E,\bnu_n)$, $n\in {\Bbb N}$, $G_{\beta}g:=\int_0^\infty e^{ 
-\beta t}T_{t}g\, dt$, $g\in L^\infty (E,\bnu)$. Since the semigroups 
$(T_{n,t})_{t\ge 0}$, $n\in {\Bbb N}$, and $(T_{t})_{t\ge 0}$ are not 
necessarily contractive, the associated families of resolvents $(G_{n, 
\beta})_{\beta >0}$, $n\in {\Bbb N}$, and $(G_\beta)_{\beta >0}$ may 
not directly be well-defined on the corresponding $L^2$-spaces. 
\medskip

Set $D(\hat{A}_n):=D(A_n)$, $n\in {\Bbb N}$, $D(\hat{A}):=D(A)$, and 
\begin{eqnarray}\label{2.14}
\hat{A}_nu_n:=A_nu_n-{\textstyle\frac12} A_n'\1\cdot u_n\, , \quad u_n\in 
D(\hat{A}_n),\ n\in {\Bbb N}, 
\end{eqnarray} 
and 
\begin{eqnarray}\label{2.15}
\hat{A}u:=Au-{\textstyle\frac12} A'\1\cdot u\, , \quad u\in D(\hat{A})\, . 
\end{eqnarray} 
Let $id_n$ denote the identity operator in $L^2(E,\bnu_n)$, $n\in {\Bbb N}$, 
and let $id$ denote the identity operator in $L^2(E,\bnu)$. 
\begin{lemma}\label{Lemma2.12} 
Suppose (c6) and let $C:={\T\frac12}\|A'\1\|_{L^\infty(E,\sbnu)}\, \vee\, 
\sup_{n\in {\Bbb N}}{\T\frac12}\|A_n'\1\|_{L^\infty (E,\sbnu_n)}$. Then we 
have the following assertions.\\ 
(a) The 
operator $\hat{A}_n$ is the generator of a strongly continuous contraction 
semigroup $(\hat{T}_{t,n})_{t\ge 0}$ in $L^2(E,\bnu_n)$, $n\in {\Bbb N}$, 
and the operator $\hat{A}$ is the generator of a strongly continuous 
contraction semigroup $(\hat{T}_t)_{t\ge 0}$ in $L^2(E,\bnu)$. In addition, 
we have $\|\hat{T}_{t,n}g\|_{L^\infty(E,\sbnu_n)}\le e^{Ct}\|g\|_{L^\infty 
(E,\sbnu_n)}$ if $g\in L^\infty(E,\bnu_n)$, $n\in {\Bbb N}$, and $\|\hat{T 
}_tg\|_{L^\infty(E,\sbnu)}\le e^{Ct}\|g\|_{L^\infty(E,\sbnu)}$ if $g\in 
L^\infty(E,\bnu)$, $t\ge 0$. \\ 
(b) Let $\alpha\ge C$. The operator $A_n-\alpha\, id_n$ is the generator of 
a strongly continuous contraction semigroup $(T_{\alpha,n,t})_{t\ge 0}$ in 
$L^2(E,\bnu_n)$, $n\in {\Bbb N}$, and the operator $A-\alpha\, id$ is the 
generator of a strongly continuous contraction semigroup $(T_{\alpha ,t} 
)_{t\ge 0}$ in $L^2(E,\bnu)$. 
\end{lemma} 
Proof. Adapting the ideas of the proof of Lemma \ref{Lemma2.11} we obtain 
for $u\in L^\infty(E,\bnu)$
\begin{eqnarray*}
&&\hspace{-.5cm}\langle u,u\rangle^{1/2}\cdot{\T\frac1t}\left(\langle T_tu, 
T_tu\rangle^{1/2}-\langle u,u\rangle^{1/2}\right) \\ 
&&\hspace{.5cm}\le\frac{\langle u,u\rangle^{1/2}}{\left(\int T'_t\1\cdot 
u^2\, d\bnu\right)^{1/2}+\langle u,u\rangle^{1/2}}\cdot{\T\frac1t}\left( 
\int T'_t\1\cdot u^2\, d\bnu-\langle u,u\rangle\right)
\end{eqnarray*} 
which results for $u=T_sv$, $v\in L^\infty(E,\bnu)$, $t\to 0$ in $\frac{d} 
{ds}\langle T_sv,T_sv\rangle^{1/2}\le\frac12\|A'\1\|\langle T_sv,T_sv 
\rangle^{1/2}$, $s\ge 0$. This and approximation of $v\in L^2(E,\bnu)$ by 
$L^\infty(E,\bnu)$-functions implies $\langle T_sv,T_sv\rangle^{1/2}\le 
e^{\frac12\|A'\1\|s}\langle v,v\rangle^{1/2}$, $v\in L^2(E,\bnu)$, $s\ge 0 
$. Similarly, $\langle T_{n,s}v,T_{n,s}v\rangle^{1/2}\le e^{\frac12\|A_n' 
\1\|s}\langle v,v\rangle^{1/2}$, $v\in L^2(E,\bnu_n)$, $s\ge 0$, $n\in 
{\Bbb N}$. 
\medskip 

\nid
(a) We prove the claim for $(\hat{T})_{t\ge 0}$. According to the 
Phillips-Lumer Theorem (cf. \cite{P83}, Theorem I.4.3, or \cite{Yo80}, 
Section IX.8) and our Lemma \ref{Lemma2.11}, it is sufficient to 
demonstrate that, for some $\hat{\alpha}>\frac12\|A'\1\|_{L^\infty(E, 
\sbnu)}$ the range of $(\hat{\alpha}\, id-\hat{A})$ is  $L^2(E,\bnu)$. 
But this follows immediately from the Feynman-Kac formula, 
\begin{eqnarray}\label{2.16}
u={\Bbb E}_\cdot\left(\int_0^\infty\exp\left\{-\hat{\alpha}t-\int_0^t 
{\textstyle\frac12}A'\1(X_s)\, ds\right\} v(X_t)\, dt\right)
\end{eqnarray} 
(${\Bbb E}$ stands for the expectation) which represents the solution to 
\begin{eqnarray*}
\hat{\alpha}u-\hat{A}u=\hat{\alpha}u -Au+{\textstyle\frac12} A'\1\cdot u 
=v\, ,\quad v\in L^2(E,\bnu)\, . 
\end{eqnarray*} 
We recall the initial step of the present proof and note that therefore 
$u=G_{\hat{\alpha}}(v-{\textstyle\frac12}A'\1\cdot u)\in D(A)$ and that 
(\ref{2.16}) corresponds to  
\begin{eqnarray}\label{2.17}
\hat{T}_tv={\Bbb E}_\cdot\left(\exp\left\{-\int_0^t{\textstyle\frac12}A' 
\1(X_s)\, ds\right\} v(X_t)\right)\, ,\quad v\in L^2(E,\bnu). 
\end{eqnarray} 
We get $\|\hat{T}_tg\|_{L^\infty(E,\sbnu)}\le e^{Ct}\|g\|_{L^\infty(E, 
\sbnu)}$ if $g\in L^\infty(E,\bnu)$, $t\ge 0$, from (\ref{2.17}) and 
condition (c6$(i),(ii)$).
\medskip 

\nid
(b) Keeping (\ref{2.16}) and (\ref{2.17}) in mind, this follows from similar 
considerations noting that, as a consequence Lemma \ref{Lemma2.11}, $S_{n, 
\alpha}$, $n\in {\Bbb N}$, and $S_\alpha$ are non-negative forms. 
\qed 
\medskip 

\nid 
{\bf Remark} (8) Another consequence of Lemma \ref{Lemma2.12} is that, 
besides the definitions (\ref{2.12})-(\ref{2.15}), $\hat{A}_n$, $\hat{S}_n$, 
$(\hat{T}_{n,t})_{t\ge 0}$, $n\in {\Bbb N}$, and $\hat{A}$, $\hat{S}$, 
$(\hat{T}_{t})_{t\ge 0}$ are related as described in Section 1. 
\medskip 

\nid 
In addition, let $(\hat{G}_{n,\beta})_{\beta\ge 0}$, denote the resolvent 
associated with $\hat{A}_n$, $\hat{S}_n$, $(\hat{T}_{n,t})_{t\ge 0}$, $n\in 
{\Bbb N}$, and let $(\hat{G}_\beta)_{\beta\ge 0}$, denote the resolvent 
associated with $\hat{A}$, $\hat{S}$, $(\hat{T}_{t})_{t\ge 0}$. 
\medskip 

\nid 
In order to handle the application in \cite{Lo14-3}, it also seems to be 
beneficial to consider the following stronger condition in place of (c3) 
and the related Lemma \ref{Lemma2.13}.  
\begin{itemize} 
\item[(c3')] 
%
\begin{itemize} 
\item[$(i)$] 
If ${\cal C}\subseteq L^\infty (E,\bnu)$ then $\{G_\beta g:g\in L^\infty 
(E,\bnu)\, , \ \beta >0\}\subseteq {\cal C}$ in the sense that for every 
$g\in L^\infty (E,\bnu)$, there is a $u\in {\cal C}$ with $G_\beta g=u$ 
${\bnu}$-a.e.; otherwise, $D(S)\subseteq {\cal C}$. 
\item[$(ii)$] 
If, for $n\in {\Bbb N}$, ${\cal C}\subseteq L^\infty (E,\bnu_n)$ then $\{ 
G_{n,\beta} g:g\in L^\infty (E,\bnu_n)\, , \ \beta >0\}\subseteq {\cal C}$ 
in the sense that for every $g\in L^\infty (E,\bnu_n)$, there is a $u\in 
{\cal C}$ with $G_{n,\beta} g=u$ ${\bnu_n}$-a.e.; otherwise, $D(S_n)\subseteq 
{\cal C}$. 
\end{itemize} 
\end{itemize} 
\begin{lemma}\label{Lemma2.13} 
Suppose (c6). (a) If (c3'(i)) then condition (c3(i)) holds for $\hat{S}$ in 
place of $S$. \\ 
(b) If (c3'(ii)) then condition (c3(ii)) holds for $\hat{S}_n$, in place of 
$S_n$, $n\in {\Bbb N}$. 
\end{lemma} 
Proof. We show only (a). {\it Step 1 } Let $\beta>C$ and $g\in {\cal C}$.  
We have 
\begin{eqnarray*}
\beta G_\beta g-AG_\beta g=g=\beta\hat{G}_\beta g-A\hat{G}_\beta g+ 
{\textstyle\frac12}A'\1\cdot \hat{G}_\beta g. 
\end{eqnarray*} 
In addition, we recall that $\hat{G}_\beta g\in L^\infty (E,\bnu)$ if $g\in 
L^\infty (E,\bnu)$ which is a consequence of Lemma \ref{Lemma2.12} (a) since 
$\beta>C$. If $g\in L^2(E,\bnu)$ then Lemma \ref{Lemma2.12} (b) guarantees 
well-definiteness of $G_\beta g$ for $\beta>C$. 

Because of $A'\1\in L^\infty (E,\bnu)$ (cf. condition (c6$(ii)$), we have 
$g-{\textstyle\frac12}A'\1\cdot\hat{G}_\beta g\in L^\infty (E,\bnu)$ if 
${\cal C}\subseteq L^\infty (E,\bnu)$, or otherwise $g-{\textstyle\frac12} 
A'\1\cdot\hat{G}_\beta g\in L^2(E,\bnu)$. The above identity implies 
\begin{eqnarray*}  
\hat{G}_\beta g=G_\beta\left(g-{\textstyle\frac12}A'\1\cdot\hat{G}_\beta g 
\right)\, .   
\end{eqnarray*} 
Condition (c3$(i)$) restricted to $\beta>C$ relative to $\hat{S}$ follows. 
\medskip 

\noindent 
{\it Step 2 } Let now $\beta\in (0,C]$ and $g\in {\cal C}$. The sum $f:= 
\sum_{k=0}^\infty(C\hat{G}_{\beta+C})^kg$ converges in $L^2(E,\bnu)$ 
according to Lemma \ref{Lemma2.12} (a) and is the unique solution to $f- 
C\hat{G}_{\beta+C}f=g$. By $\beta\hat{G}_{\beta+C}f-\hat{A}\hat{G}_{\beta 
+C}f=f-C\hat{G}_{\beta+C}f=g$ we have $\hat{G}_\beta g=\hat{G}_{\beta+C}f$ 
and as in Step 1 of the present proof
\begin{eqnarray}\label{2.18} 
\hat{G}_\beta g=\hat{G}_{\beta+C}f&&\hspace{-.5cm}=G_{\beta+C}\left(f-{ 
\textstyle\frac12}A'\1\cdot\hat{G}_{\beta+C} f\right)\nonumber \\ 
&&\hspace{-.5cm}=G_{\beta+C}\left(g+(C-{\textstyle\frac12} A'\1)\cdot 
\hat{G}_{\beta+C} f\right)\nonumber \\ 
&&\hspace{-.5cm}=G_{\beta+C}\left(g+(C-{\textstyle\frac12} A'\1)\cdot 
\hat{G}_\beta g\right)\, . 
\end{eqnarray} 
If ${\cal C}\subseteq L^\infty (E,\bnu)$ then Lemma \ref{Lemma2.12} (a) 
says that $\hat{G}_\beta g\in L^\infty (E,\bnu)$. In this case, we obtain 
$g+(C-{\textstyle\frac12} A'\1)\cdot\hat{G}_\beta g\in L^\infty (E,\bnu)$ 
and with the right-hand side of (\ref{2.18}) and the first part of 
condition (c3'$(i)$) we verify $\hat{G}_\beta g\in {\cal C}$. 

If ${\cal C}\subseteq L^\infty (E,\bnu)$ does not hold then again from  
(\ref{2.18}) and condition (c3'$(i)$) we get $\hat{G}_\beta g\in{\cal C}$. 
\qed 
\begin{theorem}\label{Theorem2.14} 
Let $\, C:={\T\frac12}\|A'\1\|_{L^\infty(E,\sbnu)}\vee\, \sup_{n\in {\Bbb 
N}}{\T\frac12}\|A_n'\1\|_{L^\infty (E,\sbnu_n)}$. Suppose (c1),(c2), and 
(c6). Furthermore, suppose (c3) for $\hat{S}_n$, $n\in {\Bbb N}$, and 
$\hat{S}$ in place of $S_n$, $n\in{\Bbb N}$, and $S$. Assume that the forms 
$\hat{S}_n$, $n\in {\Bbb N}$, converge to the form $\hat{S}$ in the sense 
of Definition \ref{Definition2.4}. \\ 
(a) Let $\alpha\ge C$. Then the forms $S_{n,\alpha}$, $n\in {\Bbb N}$, 
converge to the form $S_\alpha$ in the sense of Definition 
\ref{Definition2.4}. \\ 
(b) For all $f,g\in L^2(E,\bnu)$, all sequences $f_n\in {\cal C}$ 
$w$-converging to $f$, all sequences $g_n\in {\cal C}$ $s$-converging to 
$g$, and all $\beta>0$, we have $\hat{G}_{n,\beta}f_n\wstack{n\to\infty} 
{\lra}\hat{G}_\beta f$, $\hat{G}_{n,\beta}'f_n\wstack{n\to\infty}{\lra} 
\hat{G}_\beta 'f$ and $\hat{G}_{n,\beta}g_n\sstack{n\to\infty}{\lra}\hat 
{G}_\beta g$. Furthermore for all $g\in {\cal C}$, we have $\hat{G}_{n, 
\beta}'g\sstack{n\to\infty}{\lra}\hat{G}_\beta' g$.  \\ 
(c) Suppose (c3). For $\beta>C$, the operators $G_{n,\beta}$ and $G_\beta$ 
can be continuously extended to operators $G_{n,\beta}:L^2(E,\bnu_n)\to 
L^2(E,\bnu_n)$, $n\in {\Bbb N}$, and $G_\beta :L^2(E,\bnu)\to L^2(E,\bnu)$, 
respectively. For $f,g\in L^2(E,\bnu)$, and $f_n,g_n\in{\cal C}$, $n\in 
{\Bbb N}$, as in (b), we have $G_{n,\beta}f_n\wstack{n\to\infty}{\lra} 
G_\beta f$, $G_{n,\beta}'f_n\wstack{n\to\infty}{\lra}G_\beta'f$ and 
$G_{n,\beta}g_n\sstack{n\to\infty}{\lra}G_\beta g$, $G_{n,\beta}'g_n 
\sstack{n\to\infty}{\lra}G_\beta' g$. 
\end{theorem} 
{\bf Remark} (9) Recalling the proof of Theorem \ref{Theorem2.7}, 
especially Step 4, it will turn out that it is sufficient to require 
(c3$(i)$) as well as (c3$(ii)$) instead of (c3) (for $S_n$, $n\in {\Bbb 
N}$, and $S$) and (c3$(i)$) as well as (c3$(ii)$) instead of (c3) for 
$\hat{S}_n$, $n\in {\Bbb N}$, and $\hat{S}$ if we are not interested in 
the convergence of $G_{n,\beta}'$ or $\hat{G}_{n,\beta}'$. 
\medskip 

\nid
Proof. In order to show part (a), in Step 1-3 below, we will verify 
conditions (i)-(iii) of Definition \ref{Definition2.4} for $S_{n,\alpha}$, 
$n\in {\Bbb N}$, and $S_\alpha$. In Step 4, we will verify (b) and (c). 
\medskip 

\nid 
{\it Step 1 } Let us use the symbols $\vp_n$, $\vp$, $\psi_n$, $\psi$, 
$u_n$, $u$, as in Definition \ref{Definition2.4} with $S$ and $S_n$ replaced 
with $S_\alpha$ and $S_{n,\alpha}$. Introduce $a_n:=\alpha-{\textstyle 
\frac12}A_n'\1$, $n\in {\Bbb N}$, and $a:=\alpha-{\textstyle\frac12}A'\1$. 
Since ${\cal C}$ is dense in $L^2(E,\bnu)$ (cf. Lemma \ref{Lemma2.1} (c)), 
for given $\ve> 0$, there exists $\t \vp\in {\cal C}$ such that $\langle\vp 
-\t \vp\, , \, \vp -\t \vp\rangle <\ve$. Let such $\ve$ and $\t \vp$ be 
given. Condition (i) of Definition \ref{Definition2.4} for $S_{n,\alpha}$, 
$n\in {\Bbb N}$, and $S_\alpha$ follows from 
\begin{eqnarray*}
S_{n,\alpha}(\vp_n,\vp_n)=\hat{S}_n(\vp_n,\vp_n)+\langle a_n^{1/2}\vp_n 
\, , \, a_n^{1/2}\vp_n\rangle_n\, , \quad n\in {\Bbb N}, 
\end{eqnarray*} 
\begin{eqnarray*}
S_\alpha(\vp,\vp)=\hat{S}(\vp,\vp)+\langle a^{1/2}\vp\, , \, a^{1/2}\vp 
\rangle\, ,  
\end{eqnarray*} 
condition (i) for $\hat{S}_n$, $n\in {\Bbb N}$, and $\hat{S}$, and 
\begin{eqnarray*}
\langle a_n^{1/2}\vp_n\, , \, a_n^{1/2}\vp_n\rangle_n=\langle a_n^{1/2} 
(\vp_n -\t \vp)\, , \, a_n^{1/2}(\vp_n-\t \vp)\rangle_n+2\langle a_n\t \vp 
\, , \, \vp_n\rangle_n -\langle a_n\t \vp\, , \, \t \vp\rangle_n
\end{eqnarray*} 
together with condition (c6$(i),(ii),(iii)$) and Proposition 
\ref{Proposition2.3} (c). For the alternative in (c6$(iii)$) use 
Lemma \ref{Lemma2.1} (e) and Proposition \ref{Proposition2.3} 
(g). 
\medskip 

\nid 
{\it Step 2 } For condition (ii) of Definition \ref{Definition2.4} 
for $S_{n,\alpha}$, $n\in {\Bbb N}$, and $S_\alpha$ we take the same 
decomposition as in Step 1 for $\psi_n$ and $\psi$ instead of $\vp_n$ 
and $\vp$, $n\in {\Bbb N}$. It is then an immediate consequence of (ii) 
for $\hat{S}_n$, $n\in {\Bbb N}$, and $\hat{S}$ together with condition 
(c6$(i),(ii),(iii)$) and Proposition \ref{Proposition2.3} (c). For the 
alternative we use again Lemma \ref{Lemma2.1} (e) and Proposition 
\ref{Proposition2.3} (g). 
\medskip 

\nid 
{\it Step 3 } In order to verify condition (iii) of Definition 
\ref{Definition2.4} for $S_{n,\alpha}$, $n\in {\Bbb N}$, and $S_\alpha$, 
let $\t \beta >0$ and suppose $\sup_{n\in {\Bbb N}}\left\langle A_nu_n\, 
,\, A_nu_n\right\rangle_n<\infty$ as well as 
\begin{eqnarray*}
S_{n,\alpha+\t \beta}(u_n,\psi_n)\stack{n\to\infty}{\lra}S_{\alpha +\t 
\beta}(u,\psi)\, . 
\end{eqnarray*} 
Introducing $v_n:=\hat{G}_{n,\t \beta}(-A_nu_n+(\alpha +\t \beta)u_n)$, 
$n\in {\Bbb N}$, and $v:=\hat{G}_{\t \beta}(-Au+(\alpha +\t \beta)u)$, 
this means nothing but $\sup_{n\in {\Bbb N}}\left\langle A_nv_n\, ,\, 
A_nv_n\right\rangle_n<\infty$ as well as 
\begin{eqnarray}\label{2.19}
\hat{S}_{n,\t \beta}(v_n,\psi_n)\stack{n\to\infty}{\lra}\hat{S}_{\t \beta} 
(v,\psi)\, . 
\end{eqnarray} 
Choosing here $\psi_n=\psi$, $n\in {\Bbb N}$, and recalling Remark (2), 
we get $-\hat{A}_nv_n+\t \beta v_n\wstack{n\to\infty}{\lra}-\hat{A}v+\t 
\beta v$. Lemma \ref{Lemma2.8} implies $v_n\wstack{n\to\infty}{\lra}v$. 
Let us conclude $u_n\wstack{n\to\infty}{\lra}u$ from this. 

Let $\rho\in {\cal C}$. We have 
\begin{eqnarray}\label{2.20}
\langle u_n+\hat{G}_{n,\t \beta}(a_n\cdot u_n)\, ,\, \rho\rangle_n 
&&\hspace{-.5cm}=\langle\hat{G}_{n,\t \beta}(-A_nu_n+(\alpha +\t \beta) 
\cdot u_n)\, ,\, \rho\rangle_n \nonumber \\ 
&&\hspace{-.5cm}=\langle v_n\, ,\, \rho\rangle_n \nonumber \\ 
&&\hspace{-1cm}\stack{n\to\infty}{\lra}\langle v\, ,\, \rho\rangle 
\nonumber \\ 
&&\hspace{-.5cm}=\langle u+\hat{G}_{\t \beta}(a\cdot u)\, ,\, \rho\rangle 
\, . 
\end{eqnarray} 
On the other hand, let us assume that $u_n\wstack{n\to\infty}{\lra}\t 
u$ for some $\t u\in L^2(E,\bnu)$. From condition (c6$(i),(ii),(iii)$) 
and Proposition \ref{Proposition2.3} (c),(f) we obtain the existence of 
$\t \Psi_n\in {\cal C}$ with $\t \Psi_n=a_n\cdot u_n$ on $E_n\setminus 
N_n$ and $\t \Psi_n\wstack{n\to\infty}{\lra}a\cdot\t u$. It follows from 
Theorem \ref{Theorem2.7} (a) and the hypotheses of the present theorem 
that 
\begin{eqnarray}\label{2.21}
\langle u_n+\hat{G}_{n,\t \beta}(a_n\cdot u_n)\, ,\, \rho\rangle_n\stack 
{n\to\infty}{\lra}\langle\t u+\hat{G}_{\t \beta}(a\cdot\t u)\, ,\, \rho 
\rangle\, ,\quad \rho\in {\cal C}. 
\end{eqnarray} 
According to (\ref{2.20}) and (\ref{2.21}), it holds that $\t u-u=-\hat 
{G}_{\t \beta}(a\cdot(\t u-u))$ which implies $\t u\in D(A)$ since $u\in 
D(A)$ by hypothesis and $\hat{A}(\t u-u)-\t \beta(\t u-u)=a\cdot(\t u-u)$. 
The latter implies 
\begin{eqnarray*}
\langle-\hat{A}(\t u-u)\, ,\, \t u-u\rangle=\langle({\T\frac12}A'\1-\alpha 
-\t \beta)(\t u-u)\, ,\, \t u-u\rangle\, .  
\end{eqnarray*} 
The left-hand side is non-negative by Lemma \ref{Lemma2.11}. Since $\alpha 
\ge\frac12\|A'\1\|_{L^\infty(E,\sbnu)}$ it follows that $u=\t u$. 

Together with the above derived relation $v_n\wstack{n\to\infty}{\lra}v$ 
we now also have $u_n\wstack{n\to\infty}{\lra}u$. We recall that from 
condition (c6$(i),(ii),(iii)$) and Proposition \ref{Proposition2.3} (c), 
(f) we obtain the existence of $\bar{\Psi}_n\in {\cal C}$ with $\bar{\Psi 
}_n=a_n\cdot u_n$ on $E_n\setminus N_n$ and $\bar{\Psi}_n\wstack{n\to 
\infty}{\lra}a\cdot u$. Together with (\ref{2.19}) and Proposition 
\ref{Proposition2.3} (c), this leads to 
\begin{eqnarray*}
\hat{S}_{n,\t \beta}(u_n,\psi_n)&=&\hat{S}_{n,\t \beta}(v_n,\psi_n)-\langle 
a_n\cdot u_n\, , \, \psi_n\rangle_n \\ 
&&\hspace{-1.2cm}\stack{n\to\infty}{\lra}\hspace{-.2cm}\hat{S}_{\t \beta} 
(v,\psi)-\langle a\cdot u\, , \, \psi\rangle \\ 
&=&\hat{S}_{\t \beta}(u,\psi)\, . 
\end{eqnarray*} 
Condition (iii) of Definition \ref{Definition2.4} for $\hat{S}_n$, $n\in 
{\Bbb N}$, and $\hat{S}$ implies now that 
\begin{eqnarray}\label{2.22}
\hat{S}_{n,\t \beta}(\psi_n,u_n)\stack{n\to\infty}{\lra}\hat{S}_{\t \beta} 
(\psi,u)\, . 
\end{eqnarray} 
Since (\ref{2.22}) holds for arbitrary $D(S_n)\cap{\cal C}\ni\psi_n\sstack 
{n\to\infty}{\lra}\psi\in D(S)$, it is also true for $\psi_n$ replaced with 
$\psi_n+\hat{G}_{n,\t \beta}(\Psi_n)$ and $\psi$ replaced with $\psi+\hat{G 
}_{\t \beta}(a\cdot\psi)$ where $\Psi_n\equiv\Psi_n(\psi)\in {\cal C}$ with 
$\Psi_n=A_n'\1\cdot\psi$ on $E_n\setminus N_n$ and $E_n$ is as in (c6$(iii)$). 
For $\psi_n+\hat{G}_{n,\t \beta}(\Psi_n)\in{\cal C}$, $n\in {\Bbb N}$, recall 
conditions (c6$(iii)$) and (c3) for $\hat{S}_n$ as well as $\hat{S}$, and 
Lemma \ref{Lemma2.1} (e). For $\psi_n +\hat{G}_{n,\t \beta}(\Psi_n)\sstack{n 
\to\infty}{\lra}\psi+\hat{G}_{\t \beta}(a\cdot\psi)$, consult condition 
(c6$(i),(ii),(iii)$), Proposition \ref{Proposition2.3} (g), and Theorem 
\ref{Theorem2.7}. From (\ref{2.22}), we obtain 
\begin{eqnarray*}
S_{n,\alpha+\t \beta}(\psi_n,u_n)&=&\hat{S}_{n,\t \beta}(\psi_n+\hat{G 
}_{n,\t \beta}(a_n\cdot\psi_n),u_n) \\ 
&=&\hat{S}_{n,\t \beta}(\psi_n+\hat{G}_{n,\t \beta}(\Psi_n),u_n)+\langle 
a_n\cdot\psi_n-\Psi_n\, ,\, u_n\rangle_n \vphantom{\int^1} \\ 
&&\hspace{-1.2cm}\stack{n\to\infty}{\lra}\hspace{-.2cm}\hat{S}_{\t \beta} 
(\psi +\hat{G}_{\t \beta}(a\cdot\psi),u) \\ 
&=&S_{\alpha+\t \beta}(\psi,u)\, ; \vphantom{\int^1}
\end{eqnarray*} 
for $\langle a_n\cdot\psi_n -\Psi_n\, ,\, u_n\rangle_n\stack{n\to\infty} 
{\lra}0$ we have already mentioned that $\sup_{n\in {\Bbb N}}\langle u_n\, 
,\, u_n\rangle_n <\infty$. Thus, we have (iii) for $S_{n,\alpha}$, $n\in 
{\Bbb N}$, and $S_\alpha$. Part (a) has been verified. 
\medskip 

\nid
{\it Step 4 } Part (b) is a corollary of the Lemmas \ref{Lemma2.11} and 
\ref{Lemma2.12} (a) and Theorem \ref{Theorem2.7}. For part (c), let $\beta>C$ 
and recall Lemma \ref{Lemma2.12} (b). The claim follows from Part (a) of 
the present theorem and from Theorem \ref{Theorem2.7}. 
\qed 
\medskip 

\nid
{\bf Remark} (10) It follows from $\, C={\T\frac12}\|A'\1\|_{L^\infty(E, 
\sbnu)}\vee\, \sup_{n\in {\Bbb N}}{\T\frac12}\|A_n'\1\|_{L^\infty (E, 
\sbnu_n)}<\infty$, relation (\ref{2.17}), and the contractivity of the 
semigroups $(\hat{T}_{t,n})_{t\ge 0}$ in $L^2(E,\bnu_n)$, $n\in {\Bbb N}$, 
and $(\hat{T}_t)_{t\ge 0}$ in $L^2(E,\bnu)$, cf. Lemma \ref{Lemma2.12} (a), 
that 
\begin{eqnarray*}
\langle T_{t,n}v\, ,\, T_{t,n}v\rangle_n\le e^{2Ct}\langle v,v\rangle_n\, 
,\quad n\in {\Bbb N},\quad{\rm and}\quad\langle T_tv\, ,\, T_tv\rangle\le 
e^{2Ct}\langle v,v\rangle\, ,\quad t\ge 0.
\end{eqnarray*} 
Let us assume (c1), (c2), (c3$(i),(ii)$) and (c5). Keeping Theorem 
\ref{Theorem2.14} (c) and in mind replacing in the proof of 
Theorem \ref{Theorem2.10} and in Remark (7) contractivity by this property, 
for $\beta>C$ and ${\cal C}\ni g_n\sstack{n\to\infty}{\lra}g\in {\cal C}$ 
we obtain the following. $G_{n,\beta}g_n\sstack{n\to\infty}{\lra}G_\beta g$ 
iff $G_{n,\beta}g\sstack{n\to\infty}{\lra}G_\beta g$ iff $\ T_{n,t}g\sstack 
{n\to\infty}{\lra} T_tg$ iff $\ T_{n,t}g_n\sstack{n\to\infty}{\lra} T_tg$. 
 
\section{Relative Compactness} 
\setcounter{equation}{0}

In this section, we specify the setting of Section 2. For this, 
let ${\cal M}_1(\overline{D})$ denote the space of all probability 
measures on $(\overline{D},{\cal B}(\overline{D}))$ where $D$ is a 
bounded $d$-dimensional domain or, more general, a bounded 
$d$-dimensional Riemannian manifold for some $d\in {\Bbb N}$. In 
addition, let ${\cal M}_\partial(\overline{D})$ be the set of all 
equivalence classes $\mu$ such that $m_1,m_2\in\mu$ implies 
$m_1|_D=m_2|_D$. In this section we will assume that $E$ is one 
of the spaces ${\cal M}_1(\overline{D})$ or ${\cal M}_\partial( 
\overline{D})$. In the case of $E={\cal M}_\partial(\overline{D})$, 
we identify all points belonging to $\partial D$ with each other. 
By $r(x,y):=|x-y|\wedge\left(\inf_{b\in\partial D}|b-x|+\inf_{b\in 
\partial D}|b-y|\right)$ if $x,y\in D$ and $r(x,\partial D)=r( 
\partial D,x):=\inf_{b\in\partial D}|b-x|$ if $x\in D$, as well as 
$r(\partial D,\partial D):=0$ the space $(D\cup\partial D,r)$ 
becomes a separable, complete, and compact metric space. Furthermore, 
continuity on $D$ with respect to $r$ coincides with continuity 
with respect to the Euclidean metric and $\{f\in C(\overline{D}):f$ 
constant on $\partial D\}$ is the set of all continuous functions 
on $(D\cup\partial D,r)$. We would like to refer to similarities of 
this construction to that in \cite{Ma89}. 

Let both spaces ${\cal M}_1(\overline{D})$ and ${\cal M}_\partial 
(\overline{D})$ be endowed with the Prokhorov metric. We note that 
in this way ${\cal M}_1(\overline{D})$ and ${\cal M}_\partial 
(\overline{D})$ are separable, complete, and compact spaces.
 
Furthermore, for $n\in {\Bbb N}$, let $E_n'$ be the set of all 
measures $\mu$ in $E$ of the form $\mu=\frac{1}{n}\sum_{i=1}^n 
\delta_{z_i}$ where $z_1,\ldots ,z_n \in\overline{D}$ and $\delta_z$ 
denotes the Dirac measure concentrated at $z$. Furthermore, let $E_1 
:=E_1'$, $E_{n+1}:=E_{n+1}'\setminus\bigcup_{i=1}^nE_i$, $n\in 
{\Bbb N}$, and $E_0:=E\setminus\bigcup_{n=1}^\infty E_n$. According 
to the basic setting of Subsection 2.1 $E_n$ and $E_n'$ differ by 
$\bnu_n$-null set, $n\in {\Bbb N}$. It is therefore reasonable to 
{\it identify} $L^p(E,\bnu_n)$ {\it with both} $L^p(E_n,\bnu_n)$ and 
$L^p(E_n',\bnu_n)$, $1\le p\le\infty$, $n\in {\Bbb N}$. 

To be consistent with \cite{Lo13}, we will keep on writing $C_b(E)$
for $C(E)$. Choose ${\cal F}:=C_b(E)$ and note that therefore ${\cal 
C}$ is now the space of all functions $\vp\in {\cal D}$ satisfying the 
following.  
\begin{itemize} 
\item[(c1')] $\vp$ is bounded and continuous on $E_n$, $n\in {\Bbb N}$. 
\item[(c2')] $\langle \vp\, , \, \psi\rangle_n\stack{n\to\infty}{\lra} 
\langle\vp\, , \, \psi\rangle$ for all $\psi\in C_b(E)$.  
\end{itemize}

Obviously, ${\cal F}$ is dense in $L^2(E,\bnu)$. Under the above choice 
of ${\cal F}$, compatibility with Section 2 and well-definiteness of ${\cal 
C}$ is subject to the subsequent. 
\begin{proposition}\label{Proposition3.1} 
${\cal F}=C_b(E)$ is a subset of ${\cal D}$ if and only if the measures 
$\bnu_n$, $n\in {\Bbb N}$, are weakly convergent to $\bnu$ as $n\to\infty$; 
in symbols $\bnu_n\stack{n\to \infty}{\Ra}\bnu$. 
\end{proposition}

Let us assume that $\bnu_n\stack{n\to \infty}{\Ra}\bnu$ and note that (c1') 
and (c2') are now the defining properties of ${\cal C}\subseteq {\cal D}$. 
\medskip 

%
As in Section 2, let us assume that there are Markov processes $X^n$ and 
$X$ associated with the semigroups $(T_{n,t})_{t\ge 0}$, $n\in {\Bbb N}$, 
and $(T_t)_{t\ge 0}$. Suppose that the paths of the processes $X^n$, 
$n\in {\Bbb N}$, are cadlag. 

Define the measures $P_{{\sbnu}_n}:=\int_E P^n_\mu\, {\bnu}_n(d\mu)$, 
$n\in{\Bbb N}$, and $P_{\sbnu}:=\int_E P_\mu\, {\bnu}(d\mu)$, and 
introduce the processes ${\bf X}^n=((X^n_t)_{t\ge 0},P_{{\sbnu}_n})$ 
and ${\bf X}=((X_t)_{t\ge 0},P_{\sbnu})$. Moreover, let $\, {\Bbb E}^n_{ 
\mu}$ be the expectation corresponding to $P^n_\mu$, $\mu\in E_n$, and 
let $\, {\Bbb E}_{{\sbnu}_n}$ be the expectation corresponding to $P_{ 
{\sbnu}_n}$, $n\in {\Bbb N}$. Let us introduce the set of test functions 
we are going to work with in this section. Suppose the following.
\begin{itemize}
\item[(c7)] There exists an algebra $\t C_b(E)\subseteq C_b(E)$ of 
everywhere on $E$ defined functions with $\t C_b(E)\subseteq{\cal G}$ 
in the sense that, for every $f\in \t C_b(E)$, there is a $g\equiv g(f) 
\in {\cal C}$ and a $\beta>0$ with $f=\beta G_\beta g$ ${\bnu}$-a.e. 
$\t C_b(E)$ contains the constant functions and separates points in $E$.
\end{itemize}
{\bf Remark} (1) Note that, for $f\in \t C_b(E)$, the existence of one 
$g\equiv g(f)\in {\cal C}$ and one $\beta>0$ such that $f=\beta G_\beta g$ 
${\bnu}$-a.e. implies that, for all $\beta'>0$, there is a $g'\equiv g'(f, 
\beta')\in {\cal C}$ such that $f=\beta' G_{\beta'} g'$ ${\bnu}$-a.e. This 
follows from $Af=\beta f-\beta g\in {\cal C}$ and $g':=g+(1/\beta-1/\beta' 
)Af=(g+(Af-\beta f)/\beta)-(Af-\beta' f)/\beta'=-(Af-\beta' f)/\beta'$. 
\medskip

\nid
For $f\in\t C_b(E)$, $g=g(f)\in {\cal C}$, and a given sequence $\ve_n>0$, 
$n\in {\Bbb N}$, introduce 
\begin{eqnarray}\label{3.1}
B:=\bigcup_{n=1}^\infty\left\{\mu\in E_{n}:\left|\beta G_{n,\beta} 
g(\mu )-f(\mu )\right|\ge\ve_n\|g\|\right\}\, . 
\end{eqnarray}
Furthermore, let $\tau_{B^c}\equiv\tau_{B^c}^n(g)$ denote the first exit 
time of ${\bf X}^{n}$ from the set $B^c\cap E_{n}$, $n\in {\Bbb N}$. 
Let $T>0$ and set 
\begin{eqnarray*}
\gamma_n\equiv\gamma_n(f):=\sup_{s\in [0,T+1]}\left|\beta G_{n,\beta} 
g(X^n_s)-f(X^n_s)\right|\, , \quad n\in {\Bbb N}. 
\end{eqnarray*}
In order to prove relative compactness of the families of processes 
$f({\bf X}^n)=((f(X^n_t))_{t\ge 0},$ $P_{{\sbnu}_n}\circ f^{-1})$, 
$n\in{\Bbb N}$, we need one more technical condition. In particular, 
we specify the sequence $\ve_n>0$, $n\in {\Bbb N}$. 
\begin{itemize}
\item[(c8)] There is a sequence $\ve_n>0$, $n\in {\Bbb N}$, with $\ve_n\stack 
{n\to\infty}{\lra}0$ such that with $B\equiv B((\ve_n)_{n\in {\Bbb N}})$ 
defined in (\ref{3.1})
\begin{eqnarray*}
\, {\Bbb E}_{{\sbnu}_n}\left(e^{-\beta\tau_{B^c}}\right)\stack{n\to\infty} 
{\lra}0
\end{eqnarray*}
whenever $\langle f-\beta G_{n,\beta}g\, ,\, f-\beta G_{n,\beta}g\rangle_{n}
\stack{n\to\infty}{\lra}0$. 
\end{itemize}
\begin{theorem}\label{Theorem3.2} 
(a) Let the following be satisfied: 
\begin{itemize} 
\item[(i)] Conditions (c3), (c7), and (c8) hold. 
\item[(ii)] We have the hypotheses of Theorem \ref{Theorem2.14}, namely 
\begin{itemize} 
\item[{}] (c3) for $\hat{S}_n$, $n\in {\Bbb N}$, and $\hat{S}$ in place of 
$S_n$, $n\in {\Bbb N}$, and $S$, 
\item[{}] (c6), 
\item[{}] the forms $\hat{S}_n$, $n\in {\Bbb N}$, converge to the form 
$\hat{S}$ in the sense of Definition \ref{Definition2.4}. 
\end{itemize}
\end{itemize}
Then, for $f\in \t C_b(E)$, the family of processes $f({\bf X}^n)=((f( 
X^n_t))_{t\ge 0},P_{{\sbnu}_n}\circ f^{-1})$, $n\in{\Bbb N}$, is 
relatively compact with respect to the topology of weak convergence of 
probability measures over the Skorokhod space $D_{[-\|f\|,\|f\|]}([0, 
\infty))$. \\ 
(b) The family of processes ${\bf X}^n=((X^n_t)_{t\ge 0}$, $P_{{\sbnu}_n} 
)$, $n\in {\Bbb N}$, is relatively compact with respect to the topology 
of weak convergence of probability measures over the Skorokhod space $D_E 
([0,\infty))$. 
\end{theorem}
Proof. (a) We will apply Theorem 3.8.6 of S.N. Ethier, T. Kurtz \cite{EK86}. 
For $n\in {\Bbb N}$ and $t\ge 0$, let ${\cal F}^n_t$ denote the 
$\sigma$-algebra generated by the family $(X^n_s)_{0\le s\le t}$. In 
Steps 1 and 2 below, we will keep $n\in {\Bbb N}$ fixed. In Step 3, we 
will then pass to the limit as $n\to\infty$. 
\medskip

\nid 
{\it Step 1 } Let $f\in\t C_b(E)$ and $\beta >0$. Because of (c7) and 
Remark (1), there exist $g_1,g_2\in {\cal C}$ with $f^2=\beta G_{\beta} 
g_1$ ${\bnu}$-a.e. and $f=\beta G_{\beta} g_2$ ${\bnu}$-a.e. For $0<\delta 
<1$ and $0\le t\le T$, $0\le u\le \delta$, and $\beta >0$, we have 
\begin{eqnarray}\label{3.2}
&&\hspace{-1cm} 
\, {\Bbb E}_{{\sbnu}_n}\left[\left(f(X^{n}_{t+u})-f(X^n_t)\right)^2 
|{\cal F}^n_t\right] \nonumber \vphantom{\sum_k^2}\\ 
&&\hspace{-.5cm} =\, {\Bbb E}_{{\sbnu}_n}\left[\left. (f(X^n_{t+u}) 
)^2-(f(X^{n}_t))^2\right|{\cal F}^{n}_t\right]-2f(X^{n}_t)\, {\Bbb 
E}_{{\sbnu}_n}\left[\left. f(X^{n}_{t+u})-f(X^{n}_t)\right|{\cal 
F}^n_t\right]\nonumber\vphantom{\sum_k} \\ 
&&\hspace{-.5cm} \le\left|\, {\Bbb E}_{{\sbnu}_n}\left[\left. f^2 
(X^{n}_{t+u})-\beta G_{n,\beta}g_1(X^n_{t+u})\right|{\cal F}^{n}_t 
\right]\right|+\left|\, {\Bbb E}_{{\sbnu}_n}\left[\left.\beta G_{n, 
\beta}g_1(X^n_{t+u})\right.\right.\right. \nonumber \vphantom{\sum_k} 
\\  
&&\hspace{2cm}\left.\left.\left.-\beta G_{n,\beta}g_1(X^n_t) 
\right|{\cal F}^n_t\right]\right|+\left|\, {\Bbb E}_{{\sbnu}_n}\left[ 
\left.\beta G_{n,\beta}g_1(X^n_t)-f^2(X^n_t) \right|{\cal F}^n_t 
\right]\right| \nonumber \vphantom{\sum_k} \\  
&&\hspace{-.5cm} \hphantom{\le}+2\|f\|\left|\, {\Bbb E}_{{\sbnu}_n} 
\left[\left.f(X^{n}_{t+u})-\beta G_{n,\beta} g_2(X^n_{t+u})\right| 
{\cal F}^{n}_t\right]\right|+2\|f\|\left|\, {\Bbb E}_{{\sbnu}_n} 
\left[\left.\beta G_{n,\beta}g_2(X^n_{t+u})\right.\right.\right. 
\nonumber\vphantom{\sum_k} \\  
&&\hspace{2cm}\left.\left.\left.-\beta G_{n,\beta}g_2(X^n_t)\right| 
{\cal F}^n_t\right]\right|+2\|f\|\left|\, {\Bbb E}_{{\sbnu}_n}\left[ 
\left.\beta G_{n,\beta}g_2(X^n_t)-f(X^n_t) \right|{\cal F}^n_t\right] 
\right|\nonumber \vphantom{\sum_k} \\ 
&&\hspace{-.5cm} \le 2\, {\Bbb E}_{{\sbnu}_n}\left[\left.\sup_{s\in 
[0,T+1]}\left|\beta G_{n,\beta} g_1(X^n_s)-f^2(X^{n}_s)\right|\right| 
{\cal F}^{n}_t\right]\nonumber \\ 
&&\hspace{4.5cm}+\, {\Bbb E}_{{\sbnu}_n}\left[\sup_{r\in [0,T]}\left. 
\int_r^{r+\delta}\left|A_n(\beta G_{n,\beta}g_1)(X^{n}_s)\right| 
\, ds\right|{\cal F}^n_t\right]\nonumber \\
&&\hspace{-.5cm} \hphantom{\le}+4\|f\|\, {\Bbb E}_{{\sbnu}_n}\left[ 
\left.\sup_{s\in [0,T+1]}\left|\beta G_{n,\beta} g_2(X^n_s)-f(X^n_s) 
\right|\right|{\cal F}^n_t\right]\nonumber \\ 
&&\hspace{4.5cm}+2\|f\|\, {\Bbb E}_{{\sbnu}_n}\left[\sup_{r\in [0,T]} 
\left.\int_r^{r+\delta}\left|A_n(\beta G_{n,\beta}g_2)(X^n_s) 
\right|\, ds\right|{\cal F}^n_t\right]\nonumber \\ 
&&\hspace{-.5cm}=2\, {\Bbb E}_{{\sbnu}_n}\left[\left.\gamma_n(f^2) 
\right|{\cal F}^n_t\right]+\, {\Bbb E}_{{\sbnu}_n}\left[\sup_{r\in 
[0,T]}\left.\beta\int_r^{r+\delta}\left|g_1(X^n_s)-\beta G_{n,\beta} 
g_1(X^n_s)\right|\, ds\right|{\cal F}^n_t\right] \nonumber \\ 
&&\hphantom{\le}\hspace{-.5cm}+4\|f\|\, {\Bbb E}_{{\sbnu}_n}\left[ 
\left.\gamma_n(f)\right|{\cal F}^n_t\right]+2\|f\|\, {\Bbb E}_{{ 
\sbnu}_n}\left[\sup_{r\in [0,T]}\left.\beta\int_r^{r+\delta}\left|g_2 
(X^n_s)-\beta G_{n,\beta}g_2(X^n_s)\right|\, ds\right|{\cal F}^n_t 
\right]\hspace{-6mm}\nonumber \\ 
&&\hspace{-.5cm}\le 2\, {\Bbb E}_{{\sbnu}_n}\left[\left.\gamma_n(f^2) 
\right|{\cal F}^n_t\right]+2\beta\delta\|g_1\|+4\|f\|\, {\Bbb E}_{{ 
\sbnu}_n}\left[\left.\gamma_n(f)\right|{\cal F}^{n}_t\right]+4\beta 
\delta\|f\|\|g_2\|\vphantom{\sum_k^3}\, . 
\end{eqnarray}
{\it Step 2 }For some arbitrary sequence $\ve_n>0$, $n\in {\Bbb N}$, 
satisfying $\ve_n\stack{n\to\infty}{\lra}0$ and $\tau_{B^c}\equiv\tau_{ 
B^c}^n(g_2)$ we have 
\begin{eqnarray}\label{3.3}
\, {\Bbb E}_{{\sbnu}_n}\gamma_n(f)&=&\, {\Bbb E}_{{\sbnu}_n}\left(\chi_{\{ 
\gamma_n\le\ve_n\|g_2\|\}}\gamma_n(f)\right)+\, {\Bbb E}_{{\sbnu}_n}\left( 
\chi_{\{\gamma_n>\ve_n\|g_2\|\}}\gamma_n(f)\right)\vphantom{\sum} 
\nonumber \\ 
&\le&\ve_n\|g_2\|+(\|f\|+\|g_2\|)\, {\Bbb E}_{{\sbnu}_n}\left(\chi_{\{ 
\gamma_n>\ve_n\|g_2\|\}}\right) \vphantom{\sum}\nonumber \\ 
&=&\ve_n\|g_2\|+(\|f\|+\|g_2\|)e^{\beta(T+1)}\, {\Bbb E}_{{\sbnu}_n}\left( 
e^{-\beta(T+1)}\chi_{\{\gamma_n>\ve_n\|g_2\|\}}\right)\vphantom{\sum^1} 
\nonumber \\ 
&\le&\ve_n\|g_2\|+(\|f\|+\|g_2\|)e^{\beta(T+1)}\, {\Bbb E}_{{\sbnu}_n}\left( 
e^{-\beta\tau_{B^c}}\chi_{\{\gamma_n>\ve_n\|g_2\|\}}\right)\, .\vphantom{ 
\sum^1_1}
\end{eqnarray}
%
It follows now from $f=\beta G_\beta g_2$ ${\bnu}$-a.e. and $f-\beta 
G_\beta g_2\in {\cal C}$ (cf. (c3) and (c7)) that $\langle f-\beta G_{ 
\beta}g_2\, , \, f-\beta G_{\beta}g_2\rangle_{n}\stack{n\to\infty}{\lra}0$. 
Together with $G_{n,\beta}g_2\sstack{n\to\infty}{\lra}G_{\beta}g_2$ (cf. 
Theorem \ref{Theorem2.14} (c)), this leads to 
\begin{eqnarray*}
\langle f-\beta G_{n,\beta}g_2\, , \, f-\beta G_{n,\beta}g_2\rangle_{n}
\stack{n\to\infty}{\lra}0\, . 
\end{eqnarray*}
Relation (\ref{3.3}) and condition (c8) imply now 
\begin{eqnarray}\label{3.4}
\, {\Bbb E}_{{\sbnu}_n}\gamma_{n}(f)\stack{n\to\infty}{\lra}0
\quad \mbox{\rm and similarly}\quad \, {\Bbb E}_{{\sbnu}_n}\gamma_{n}(f^2) 
\stack{n\to\infty}{\lra}0\, . 
\end{eqnarray}
{\it Step 3 } Setting 
\begin{eqnarray*}
\g_n(\delta):=2\gamma_n(f^2)+4\|f\|\gamma_n(f)+2\beta\delta\|g_1\|+
4\beta\delta\|f\|\|g_2\|\, , \quad n\in {\Bbb N}, 
\end{eqnarray*}
and taking into consideration (\ref{3.2}), we observe that  
\begin{eqnarray*}
\, {\Bbb E}_{{\sbnu}_n}\left[\left(f(X^{n}_{t+u})-f(X^{n}_t)\right)^2 
|{\cal F}^{n}_t\right]\le \, {\Bbb E}_{{\sbnu}_n}\left[\left. \g_{n} 
(\delta)\right|{\cal F}^{n}_t\right]\, , \quad f\in \t C_b(E), \ 
n\in {\Bbb N}, 
\end{eqnarray*}
and from (\ref{3.4}), we obtain 
\begin{eqnarray*}
\lim_{\delta\to 0}\lim_{n\to 0}\, {\Bbb E}_{{\sbnu}_n}(\delta)=0\, .
\end{eqnarray*}
Relative compactness of the family $f({\bf X}^n)$, $n\in {\Bbb N}$, 
follows now from \cite{EK86}, Theorems 3.7.2 and 3.8.6, Remark 3.8.7, 
and the fact that the processes $f({\bf X}^n)$, $n\in {\Bbb N}$, take 
values in the compact interval $[\inf f,\sup f]$. \\ 
(b) Let us recall that $E$ is compact.  According to condition (c7), the 
Stone-Weierstrass Theorem, implies that $\t C_b(E)$ is dense in $C_b(E)$ 
with respect to the $\sup$-norm. Now the claim is a consequence of part 
(a) and \cite{EK86}, Theorem 3.9.1.
\qed


\small

\begin{thebibliography}{aa}
\addcontentsline{toc}{section}{\protect\numberline{6}{References}} 

\bibitem{CGLW12}{\sc F. Conrad, M. Grothaus, J. Lierl, O. Wittich},
Convergence of Brownian motion with a scaled Dirac delta potential. 
{\it Proc. Edinb. Math. Soc. (Series 2)} {\bf 55} No. 2 (2012), 
403-427. 

\bibitem{EK86}{\sc S. N. Ethier, T. Kurtz}, {\it Markov processes, 
Characterization and convergence}, New York Chichester Brisbane 
Toronto Singapore: John Wiley 1986. 

\bibitem{GKR07}{\sc M. Grothaus, Y. G. Kondratiev, M. R\"ockner}
$N/V$-limit for stochastic dynamics in continuous particle systems. 
{\it Probab. Theory Related Fields} {\bf 137} No. 1-2 (2007), 
121-160. 
 
\bibitem{Hi98}{\sc M. Hino}, Convergence of non-symmetric forms.  
{\it J. Math. Kyoto Univ.} {\bf 38} No. 2 (1998), 329-341. 

\bibitem{Hi09}{\sc M. Hinz} Approximation of jump processes on 
fractals. {\it Osaka J. Math.} {\bf 46} No. 1 (2009), 141-171.

\bibitem{Ki06}{\sc P. Kim} Weak convergence of censored and 
reflected stable processes. {\it Stochastic Process. Appl.} {\bf 
116} No. 12 (2006), 1792-1814.

\bibitem{Ko05}{\sc A. V. Kolesnikov}, Convergence of Dirichlet 
forms with changing speed measures on ${\Bbb R}^d$. {\it Forum 
Math.} {\bf 17} No. 2 (2005), 225-259.

\bibitem{Ko06}{\sc A. V. Kolesnikov}, Mosco convergence of Dirichlet 
forms in infinite dimensions with changing reference measures. {\it 
J. Funct. Anal.} {\bf 230} No. 2 (2006), 382-418. 

\bibitem{Ko08}{\sc A. V. Kolesnikov}, Weak convergence of diffusion 
processes on Wiener space. {\it Probab. Theory Related Fields} {\bf 
140} No. 1-2 (2008), 1-17. 

\bibitem{KU97}{\sc K. Kuwae, T. Uemura} Weak convergence of symmetric 
diffusion processes. {\it Probab. Theory Related Fields} {\bf 109} 
No. 2 (1997), 159-182. 

\bibitem{KU96}{\sc K. Kuwae, T. Uemura} Weak convergence of symmetric 
diffusion processes II. {\it Probability theory and mathematical 
statistics (Tokyo, 1995)}, River Edge, NJ: World Sci. Publ., (1996), 
266-275.

\bibitem{KS03}{\sc K. Kuwae, T. Shioya}, Convergence of spectral 
structures: a functional analytic theory and its applcations to 
spectral geometry. {\it Comm. Anal. Geom.} {\bf 11} No. 4 (2003), 
599-673. 

\bibitem{Lo09}{\sc J.-U. L\"obus}, A stationary Fleming-Viot type 
Brownian particle system. {\it Math. Z.} {\bf 263} No. 3 (2009), 
541-581. 

\bibitem{Lo13}{\sc J.-U. L\"obus}, Weak convergence of $n$-particle 
systems using bilinear forms. {\it Milan J. Math.} {\bf 81} No. 1 
(2013), 37-77. 

\bibitem{Lo14-3}{\sc J.-U. L\"obus}, Mosco type convergence and weak 
convergence for a Fleming-Viot type particle system, Preprint (2013), 
{\tt http://www.mai.liu.se/$\sim$julob/Lo14-3.pdf} 

\bibitem{Ma89}{\sc Mandelkern}, Metrization of the one-point 
compactification. {\it Proc. Amer. Math. Soc.} {\bf 107} No. 4 (1989), 
1111-1115.

\bibitem{MR92}{\sc Z.-M. Ma and M. R\"ockner}, {\it Introduction to 
the Theory of (Non-symmetric) Dirichlet Forms}, Berlin: Springer 1992.

\bibitem{Mo94}{\sc U. Mosco}, Composite media and asymptotic Dirichlet 
forms. {\it J. Funct. Anal.} {\bf 123} No. 2 (1994), 368-421. 

\bibitem{P83}{\sc A. Pazy}, {\it  Semigroups of linear operators and 
applications to partial differential equations}, New York: Springer 
1983. 

\bibitem{Pu10}{\sc O. V. Pugachev}, On the closability and convergence 
of Dirichlet forms. {\it Proc. Steklov Inst. Math.}, {\bf 270} (2010), 
216-221.

\bibitem{Su98}{\sc W. Sun}, Weak convergence of Dirichlet processes. 
{\it Sci. China Ser. A} {\bf 41} No. 1 (1998), 8-21. 

\bibitem{St99}{\sc W. Stannat}, The theory of generalized Dirichlet forms 
and its applications in analysis and stochastics. {\it Mem. Amer. Math. 
Soc.} {\bf 142} No. 678 (1999). 

\bibitem{To06}{\sc J. M. T\"olle} Convergence of non-symmetric forms 
with changing reference measures, Diploma thesis, Faculty of Mathematics, 
University of Bielefeld, 2006. 

\bibitem{Yo80}{\sc K. Yosida}, {\it Functional Analysis}, Berlin 
Heidelberg New York: Springer 1980. 

\end{thebibliography}
\end{document}